\documentclass[journal]{IEEEtran}

\newcommand{\beq}{\begin{equation}}
\newcommand{\eeq}{\end{equation}}
\newcommand{\bsq}{\begin{subequations}}
\newcommand{\esq}{\end{subequations}}
\newcommand{\bq}{\begin{eqnarray}}
\newcommand{\eq}{\end{eqnarray}}
\newcommand{\bqn}{\begin{eqnarray*}}
\newcommand{\eqn}{\end{eqnarray*}}
\usepackage{bbm}
\usepackage[pdftex]{graphicx}
\usepackage{enumerate}
\usepackage{enumitem} 

\usepackage{mathptmx}       
\DeclareMathAlphabet{\mathcal}{OMS}{cmsy}{m}{n}
\usepackage{helvet}         
\usepackage{courier}        
\usepackage{type1cm}        
\usepackage{makeidx}         
\usepackage{graphicx}        
\usepackage{url} 
\usepackage{multicol}        
\usepackage[bottom]{footmisc}
\usepackage{ifthen}
\usepackage{subfigure}
\usepackage[usenames,dvipsnames]{color}

\makeindex             

\usepackage{graphics} 
\usepackage{graphicx} 
\usepackage{epsfig} 
\usepackage{epstopdf}
\usepackage{mathptmx} 
\usepackage{times} 
\usepackage[cmex10]{amsmath}

\usepackage{mathtools}  
\usepackage{amssymb}  
\usepackage{amsthm}

\usepackage{amsfonts}
\usepackage{cite}
\usepackage{verbatim}
\usepackage{comment}
\usepackage{algorithmic}
\usepackage{multirow}
\usepackage{cite}
\usepackage{stfloats}
\usepackage{supertabular}
\usepackage{longtable}

\DeclareGraphicsExtensions{.pdf,.png,.jpg,.eps,.ps}
\renewcommand{\arraystretch}{1.2}
\usepackage{arydshln}
\usepackage{multicol}
\usepackage{colortbl}
\theoremstyle{definition}

\newtheorem{proposition}{Proposition}

\theoremstyle{definition}
\newtheorem{definition}{Definition}

\ifCLASSINFOpdf
\else
\fi

\usepackage{comment}
\usepackage{ifthen}
\newboolean{showcomments}
\setboolean{showcomments}{true}
\newcommand{\ychen}[1]{\ifthenelse{\boolean{showcomments}}
        { \textcolor{red}{YC: #1}}}
\newcommand{\tongxin}[1]{\ifthenelse{\boolean{showcomments}}
        { \textcolor{blue}{(#1)}}{}}

\usepackage{amsmath}
\allowdisplaybreaks[4]

\usepackage[ruled]{algorithm2e}

\begin{document}
        
        %
\title{An Energy Sharing Mechanism Achieving the Same Flexibility as Centralized Dispatch}
        %
        %
        %
        
\author{Yue Chen,
        Wei Wei,
        Han Wang,
        Quan Zhou,
        and Jo{\~a}o P. S. Catal{\~a}o}
        \markboth{Journal of \LaTeX\ Class Files,~Vol.~XX, No.~X, Feb.~2019}%
        {Shell \MakeLowercase{\textit{et al.}}: Bare Demo of IEEEtran.cls for IEEE Journals}
        %



        \maketitle
        
        \begin{abstract}
Deploying distributed renewable energy at the demand side is an important measure to implement a sustainable society. However, the massive small solar and wind generation units are beyond the control of a central operator. To encourage users to participate in energy management and reduce the dependence on dispatchable resources, a peer-to-peer energy sharing scheme is proposed which releases the flexibility at the demand side. Every user makes decision individually considering only local constraints; the microgrid operator announces the sharing prices subjective to the coupling constraints without knowing users' local constraints. This can help protect privacy. We prove that the proposed mechanism can achieve the same disutility and flexibility as centralized dispatch, and develop an effective modified best-response based algorithm for reaching the market equilibrium. The concept of ``absorbable region'' is presented to measure the operating flexibility under the proposed energy sharing mechanism. A linear programming based polyhedral projection algorithm is developed to compute that region. Case studies validate the theoretical results and show that the proposed method is scalable.
        \end{abstract}
        

        %
        \IEEEpeerreviewmaketitle

\section*{Nomenclature}
\addcontentsline{toc}{section}{Nomenclature}
\subsection{Indices, Sets, and Functions}
\begin{IEEEdescription}[\IEEEusemathlabelsep\IEEEsetlabelwidth{${\underline P _{mn}}$,${\overline P _{mn}}$}]
\item[$i \in \mathcal{I}$]   Consumer $i$ in set $\mathcal{I}$.
\item[$j \in \mathcal{J}$]   Prosumer $j$ in set $\mathcal{J}$.
\item[$l \in \mathcal{L}$]   Line $l$ in set $\mathcal{L}$.
\item[$f_k(.)$] Disutility function of user $k \in \mathcal{I} \cup \mathcal{J}$.
\item[$\hat{D}_k$] Local constraint for user $k \in \mathcal{I} \cup \mathcal{J}$.
\item[$\tilde{\mathcal{P}}$] Coupling constraint for all users.
\item[$W_c^D$] Absorbable region under centralized dispatch.
\item[$W_s^D$] Absorbable region under energy sharing.
\end{IEEEdescription}

\subsection{Parameters}
\begin{IEEEdescription}[\IEEEusemathlabelsep\IEEEsetlabelwidth{${\underline P _{mn}}$,${\overline P _{mn}}$}]
\item[$I$] Number of consumers.
\item[$J$] Number of prosumers.
\item[$d_i^f$] Fixed demand of consumer $i \in \mathcal{I}$.
\item[$d_j^f$] Fixed demand of prosumer $j \in \mathcal{J}$.
\item[$a$] Sharing market sensitivity.
\item[$\pi_{kl}^w, \pi_{kl}^d$] Line flow distribution factors.
\item[$F_l$] Power flow limit of line $l \in \mathcal{L}$.
\end{IEEEdescription}

\subsection{Decision Variables}
\begin{IEEEdescription}[\IEEEusemathlabelsep\IEEEsetlabelwidth{${\underline P _{mn}}$,${\overline P _{mn}}$}]
\item[$d_i$] Elastic demand of consumer $i \in \mathcal{I}$.
\item[$d_j$] Elastic demand of prosumer $j \in \mathcal{J}$.
\item[$w_j$] Renewable output of prosumer $j \in \mathcal{J}$.
\item[$p_k^{out}$] Net demand of user $k \in \mathcal{I} \cup \mathcal{J}$.
\item[$q_k$] Sharing amount of user $k \in \mathcal{I} \cup \mathcal{J}$.
\item[$\lambda_k$] Sharing price of user $k \in \mathcal{I} \cup \mathcal{J}$.
\item[$b_k$] Bid of user $k \in \mathcal{I} \cup \mathcal{J}$.
\end{IEEEdescription}

\section{Introduction}
        %
        %
        %
        %
\IEEEPARstart{T}{he} penetration of distributed renewable generation has been growing rapidly in recent decades \cite{rahbar2016energy}, which helps pave the way towards a green and sustainable smart grid. However, the inherent volatility and intermittency of wind and solar power also threaten the security of power system energy management \cite{eltigani2015challenges}. How to maintain system reliability under a high share of renewable energy becomes a crucial topic. Related research can be roughly classified into two categories: one aims to come up with an optimal and reliable operating strategy given the uncertainty model of renewable output, such as samples \cite{shuai2018stochastic} or an uncertainty set \cite{lei2019robust}. Typical techniques include stochastic optimization \cite{shuai2018stochastic,xu2019data}, robust optimization \cite{moretti2020efficient,lei2019robust}, and distributionally robust optimization \cite{hajebrahimi2020adaptive,chen2016distributionally}. The other endeavors to quantify the system's ability to accommodate renewable energy, by the region-based geometric method \cite{zhao2014variable} or various metrics \cite{zhao2015unified}.
        
For the region-based approaches, a pioneering one is the do-not-exceed limit (DNEL) proposed in \cite{zhao2014variable}. The DNEL refers to the minimum and maximum output levels of renewable energy within which the system constraints can be satisfied. The DNEL given by an optimization model maximizing the weighted total renewable output variation range is then embedded into a flexible dispatch model. Three alternative approaches were developed to derive dispatch strategies. Extensive works have been carried out to further improve the DNEL based dispatch method by taking into account historical data \cite{li2016data,qiu2016data}, corrective topology control \cite{korad2015zonal,korad2015enhancement}, and different kinds of supporting sets \cite{ma2019distributionally}. Another well-known concept is the dispatchable region \cite{wei2014dispatchable}, which contains a set of renewable output scenarios under which the DC power flow model has at least one feasible solution. It was further extended to AC power flow models \cite{liu2019real}, and multi-energy systems \cite{chen2018convex}. 
                
The above works provide profound techniques for evaluating the system's flexibility under a \emph{centralized} manner. It means all generators and users are controlled by a central operator. However, with the prevalence of distributed energy resources, the number of participants increases dramatically and their locations become more dispersed. In such circumstances, the centralized approach is difficult to implement. For one reason, massive participants add up to the difficulty of data acquisition and privacy protection \cite{chen2018optimal}; for example, the cost coefficients, preference, and capacity limits are known only to each user but not the operator. For another, the consumption of each user is beyond the control of a central authority.        

To overcome the obstacle of centralized operation, various market mechanisms were designed to manage massive participants in a distributed manner \cite{wang2021distributed,chen2018analyzing}. Coalitional game-based approaches were used for coordinating distributed storages \cite{chakraborty2018sharing} and microgrids \cite{mei2019coalitional}. A trading mechanism with Shapley value was presented in \cite{long2019game} and compared with three classic mechanisms, i.e., bill sharing, mid-market rate, and supply-demand ratio. This Shapley value can be estimated by a coalitional stratified random sampling method \cite{han2019estimation}. The scalability of prosumers' cooperative game was improved by K-means clustering \cite{han2019improving}. A Nash bargaining model was adopted in \cite{dutta2014game} to address the charge sharing among electric vehicles. For non-cooperative game-based approaches, distributed peer-to-peer energy exchange was modeled as a generalized Nash game in \cite{le2020peer} whose set of variational equilibria coincides with the set of social optimum. A generalized demand function based energy sharing mechanism was proposed in \cite{chen2019energy} with proofs of several properties of its equilibrium. A practical bidding algorithm was given in \cite{chen2020approaching}. A decentralized algorithm was proposed in \cite{chen2020decentralized} to provide renewable predictions to consumers in a virtual power plant (VPP) with the condition of its convergence. The above works investigate the strategic behaviors of participants and the market equilibrium, but the flexibility of a certain market mechanism in accommodating renewable energy has not been well explored yet.
        
To this end, this paper proposes an energy sharing scheme and quantifies its flexibility. The contributions are two-fold:

1) \textbf{Mechanism Design}. A mechanism that coordinates the peer-to-peer energy sharing among massive users is presented. Each user makes a decision subject to its local constraints, and the operator decides the energy sharing prices solely according to the system-wide coupling constraint without users' private data. The sharing market can be described by a generalized Nash game. The existence of a generalized Nash equilibrium (GNE) which can achieve social optimum is proved. We also develop a modified best-response based algorithm for reaching the sharing market equilibrium with proof of its convergence.
        
2) \textbf{Flexibility Characterization}. We generalize the concept of ``dispatchable region'' in \cite{wei2014dispatchable} under centralized dispatch to the ``absorbable region'' under any given scheme to characterize the flexibility of the proposed energy sharing mechanism geometrically. We prove that the absorbable region under energy sharing is exactly the same as that under centralized dispatch, meaning that the system’s flexibility is retained. A linear programming based polyhedral projection algorithm is developed to generate the absorbable region.

The rest of this paper is organized as follows. A peer-to-peer energy sharing mechanism is proposed in Section II and the existence and efficiency of its equilibrium are proved; A modified best-response based algorithm is developed to reach the sharing market equilibrium. To quantify the flexibility of energy sharing, the concept of ``dispatchable region'' under centralized operation is extended to the ``absorbable region'' in Section III, and a linear programming based polyhedron projection algorithm is established to compute that region. Case studies in Section IV validate the theoretical outcomes. Conclusions are drawn in Section V.

\section{Energy Sharing Mechanism}

To deal with the practical issues in centralized operation, an energy sharing mechanism for managing massive users is developed in this section. The users in the sharing market play a generalized Nash game. We prove that the market equilibrium always exists and is socially optimal. A modified best-response based algorithm is presented to reach the market equilibrium, and its convergence condition is given.

\subsection{Problem description}

We consider the optimal dispatch of a group of massive users, including consumers indexed by $i \in \mathcal{I}:=\{1,2,...,I\}$ and prosumers indexed by $j \in \mathcal{J}:=\{1,2,...,J\}$, in a stand-alone microgrid. Consumer $i$'s fixed demand is $d_i^f$ and its elastic demand is $d_i$. Each prosumer $j \in \mathcal{J}$ is equipped with renewable units whose total output is $w_j$. Its fixed demand is $d_j^f$ and elastic demand is $d_j$. In this paper, we consider a specific type of demand response, where users can adjust their demand to minimize their disutility as in \cite{li2011optimal,samadi2012advanced}. Under a centralized manner, the microgrid operator solves problem \eqref{eq:central} to maintain power balancing with the lowest total user disutility.
\bsq
\label{eq:central}
\begin{align}
    \mathop{\min}_{d_k,\forall k \in \mathcal{I} \cup \mathcal{J}} ~ & \sum \nolimits_{k \in \mathcal{I} \cup \mathcal{J}} f_k(d_k) \label{eq:central.1} \\
    \mbox{s.t.}~& \sum \nolimits_{k \in \mathcal{I} \cup \mathcal{J}} p_k^{out}=0 \label{eq:central.2}\\
    ~ & p_k^{out} = \left\{ \begin{aligned}
    & d_k^f+d_k,\forall k \in \mathcal{I} \\
    & d_k^f+d_k-w_k ,\forall k \in \mathcal{J}
    \end{aligned}\right.:\lambda_k\label{eq:central.4}\\
    ~ & d_k \in \hat{\mathcal{D}}_k,\forall k \in \mathcal{I} \cup \mathcal{J} \label{eq:central.3}\\
    ~ & p^{out} \in \tilde{\mathcal{P}} \label{eq:central.5}
\end{align}
\esq
Here, $f_k(d_k)$ characterizes the disutility of user $k \in \mathcal{I} \cup \mathcal{J}$ caused by adjusting its elastic demand, which is a strictly convex and twice differentiable function. Constraints \eqref{eq:central.2} and \eqref{eq:central.4} represent the power balancing between the total output of renewable units and the total demand. Constraint \eqref{eq:central.3} is the local feasible set of each user $k \in \mathcal{I} \cup \mathcal{J}$, e.g., the range limit for each responsive load, represented by a closed convex set $\hat{\mathcal{D}}_k$. Constraint \eqref{eq:central.5} collects the global coupling constraints, e.g., the network power flow limit, for the net demand (defined in \eqref{eq:central.4}), represented by a closed convex set $\tilde{\mathcal{P}}$. 

Throughout the paper, we assume
	
\noindent \textbf{A1}. $\{d: \mbox{s.t.} \; \eqref{eq:central.2}-\eqref{eq:central.5} \; \mbox{are satisfied.}\} \ne \emptyset$.

The centralized dispatch performs well in some rural or isolated microgrids, where the operator can get all the information needed for operation. Nevertheless, with the growing scale of microgrids and power sector decentralization, the marketization of microgrid has become a prevalent trend \cite{aram2017microgrid}. Under this circumstance, how to protect users’ \emph{information privacy} and reduce the possible market power exploitation due to \emph{information asymmetry} is getting increasingly important. To be specific, in problem \eqref{eq:central}, we assume that the operator knows all related constraints, but in practice, the $\hat{\mathcal{D}}_k$, the disutility function $f_k(.)$, and fixed demand $d_k^f$ are usually private information only available to user $k \in \mathcal{I} \cup \mathcal{J}$. Asking for such information may jeopardize users' privacy. The user may even have the incentive to misrepresent this information to lower its disutility. Besides, each user $k$ may not have access to the coupling constraint $\tilde{\mathcal{P}}$, which is known to the operator only. Thus, the centralized model \eqref{eq:central} will encounter some challenges for microgrid management in practice, and a new approach that can perform with local information and reduce the impact of information asymmetry is desired.

\textbf{Remark:} We further distinguish three similar expressions related to the information structure, i.e. information privacy, information asymmetry, and information scarcity, for a better understanding. By saying \emph{information privacy}, we mean that some information is available to a specific party while other parties cannot get access to it. The concept of \emph{information asymmetry} stems from contract theory and economics. Though it also refers to the case where different parties have different information, it highlights more on the situation that the party with more information may misrepresent their private information to gain more profit, resulting in an imbalance market power or even a market failure \cite{laffont2009theory}. For example, in order to solve problem (1), the central coordinator needs to collect information about $f_k(.)$ for all $k \in \mathcal{I} \cup \mathcal{J}$, which are private to users. Therefore, the user may deliberately misreport these data so as to lower its disutility. An example of how this asymmetric information structure may influence the outcomes of the system is given in Section IV-A. \emph{Information scarcity} indicates that a party does not know other parties’ strategies, the only information it has access to is the outcome \cite{nioche2019coordination}. In this paper, we assume that all users and the microgrid operator know others’ strategies when making their own decisions.
\subsection{Mechanism design}

\begin{figure}[t]
        \centering
        \includegraphics[width=1.0\columnwidth]{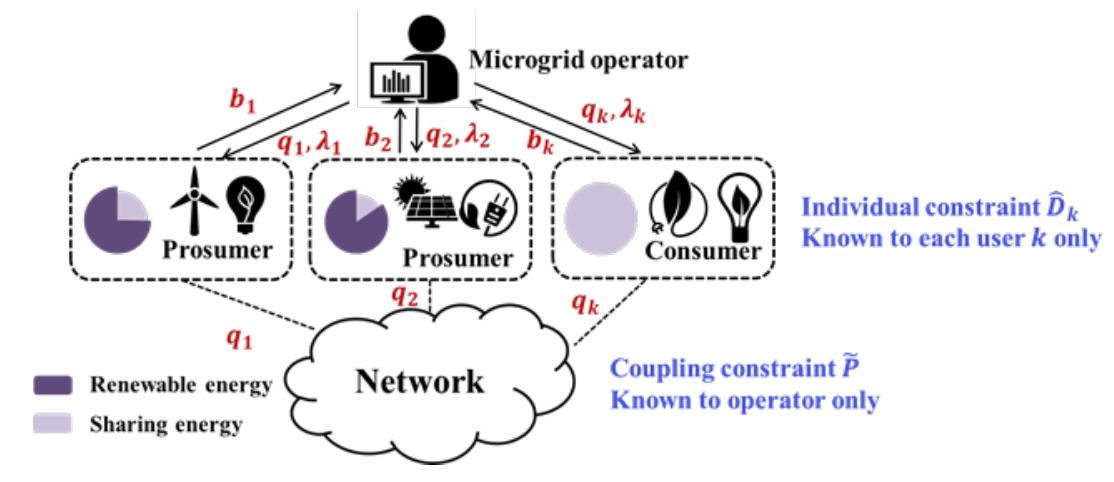}
        \caption{Sketch of the energy sharing market.}
        \label{fig:sketch}
\end{figure}
To solve problem \eqref{eq:central}, the microgrid operator needs to know $\hat D_k$, $f_k(d_k)$ and $d_k^f$, which could vary among different end-users \cite{yang2015economical}.
Therefore, the centralized operation may be impracticable due to the privacy protection requirement and possible speculative behavior under information asymmetry.
To overcome these problems, a peer-to-peer energy sharing mechanism is proposed: Each user $k \in \mathcal{I} \cup \mathcal{J}$ determines its own demand $d_k$, and meanwhile, can share energy $q_k$ with other prosumers at a price $\lambda_k$ to maintain self-power balancing. Instead of deciding on the $\lambda_k$ or $q_k$ directly, each user offers a bid $b_k$ to the operator considering its local constraint $\hat D_k$. With all $b_k,\forall k \in \mathcal{I} \cup \mathcal{J}$, the microgrid operator clears the sharing market subject to the coupling constraint $\tilde{\mathcal{P}}$ and determines the price $\lambda_k,\forall k \in \mathcal{I} \cup \mathcal{J}$; the transactive energy  $q_k:=-a\lambda_k+b_k$ is transferred via the power network, where $a>0$ represents the sharing market sensitivity. When $q_k>0$ user $k$ buys energy from the market, and when $q_k<0$ it sells. A sketch of the energy sharing market is provided in Fig.\ref{fig:sketch}. 

Under the above settings, the mathematical model of the energy sharing scheme is presented. The optimization problem of each user $k \in \mathcal{I} \cup \mathcal{J}$ is: 
\bsq
\label{eq:sharing-pro}
\begin{align}
    \mathop{\min}_{d_k,b_k}~ & f_k(d_k)+\lambda_k(-a\lambda_k+b_k) \label{eq:sharing-pro.1}\\
    \mbox{s.t.} ~& \left\{\begin{aligned}
     & -a\lambda_k+b_k=d_k^f+d_k,~\forall k \in \mathcal{I} \\
     & w_k-a\lambda_k+b_k=d_k^f+d_k,~\forall k \in \mathcal{J}
    \end{aligned}
    \right.
    \label{eq:sharing-pro.2}\\
    ~ & d_k \in \hat{\mathcal{D}}_k \label{eq:sharing-pro.3}
\end{align}
\esq
where $f_k(d_k)$ is its disutility in monetary terms, and $\lambda_k(-a\lambda_k+b_k)$ is its payment to the energy sharing market (or revenue when this term is negative). Constraint \eqref{eq:sharing-pro.2} is the power balancing condition, and \eqref{eq:sharing-pro.3} is the local limit. Note that when making a decision, each user $k \in \mathcal{I} \cup \mathcal{J}$ needs not know the coupling constraint $\tilde{\mathcal{P}}$. 
We will prove latter in Proposition \ref{prop-1} that the energy sharing prices $\lambda_k,\forall k \in \mathcal{I} \cup \mathcal{J}$ equal to the dual variables of \eqref{eq:central.2} at the optimal point, which is also the value of locational marginal prices \cite{schweppe2013spot}. When there are massive users, if we change the $b_k$ of one user, the $\lambda_k,\forall k \in \mathcal{I} \cup \mathcal{J}$  obtained by solving \eqref{eq:sharing-oper} will remain nearly unchanged. Therefore, the impact of each $b_k$ on $\lambda_k$ is negligible, and in problem \eqref{eq:sharing-pro} $\lambda_k$ is seen as an exogenously given constant. This is a common assumption in economics \cite{landsburg2013price}.

Given all users' bids $b_k,\forall k \in \mathcal{I} \cup \mathcal{J}$, the microgrid operator solves the following problem to clear the sharing market and determine the energy sharing prices.
\bsq
\label{eq:sharing-oper}
\begin{align}
    \mathop{\min}_{\lambda_k,\forall k \in \mathcal{I} \cup \mathcal{J}}~ & \sum \limits_{k \in \mathcal{I} \cup \mathcal{J}} \lambda_k^2 \label{eq:sharing-oper.1}\\
    \mbox{s.t.}~& \sum \limits_{k \in \mathcal{I} \cup \mathcal{J}} (-a\lambda_k+b_k)=0 \label{eq:sharing-oper.2}\\
    ~ & (-a\lambda+b) \in \tilde{\mathcal{P}} \label{eq:sharing-oper.3}
\end{align}
\esq
First, the sharing market needs to be cleared, which means the energy bought equals the energy sold. Therefore, we have $\sum \nolimits_{k \in \mathcal{I} \cup \mathcal{J}} q_k=0$ in condition \eqref{eq:sharing-oper.2}. Secondly, the transactive energy $q_k,\forall k \in \mathcal{I} \cup \mathcal{J}$  will be transferred via network, and should meet the network constraints \eqref{eq:sharing-oper.3}. The objective function aims to enhance market fairness by minimizing the variance of energy sharing prices, i.e.
\begin{align}
\label{eq:R1}
	~ & \sum \limits_{k \in \mathcal{I} \cup \mathcal{J}} \left(\lambda_k-\frac{1}{I}\sum \limits_{k \in \mathcal{I} \cup \mathcal{J}} \lambda_k\right)^2 = \sum \limits_{k \in \mathcal{I} \cup \mathcal{J}} \left(\lambda_k-\frac{1}{aI}\sum \limits_{k \in \mathcal{I} \cup \mathcal{J}} b_k\right)^2 \nonumber\\
	=~ & \sum \limits_{k \in \mathcal{I} \cup \mathcal{J}} \lambda_k^2 -\frac{1}{a^2I} \left(\sum \limits_{k \in \mathcal{I} \cup \mathcal{J}} b_k\right)^2
\end{align}
which is equivalent to minimizing \eqref{eq:sharing-oper.1} since the second term in \eqref{eq:R1} is a constant. Moreover, if all network constraints in \eqref{eq:sharing-oper.3} are not binding, it degenerates to the case without network constraints in \cite{chen2020approaching} with a unified price $\lambda$ given by $\lambda=\sum \nolimits_{k \in \mathcal{I} \cup \mathcal{J}} b_k /(aI)$. After clearing the market, the microgrid operator returns the price $\lambda_k$ back to each user; then the user will adjust its bid according to \eqref{eq:sharing-pro} and submit it to the operator again. This happens until convergence. We will prove later in Proposition \ref{prop-1} that the market outcome at equilibrium satisfies all local constraints $\hat{\mathcal{D}}_k$ spontaneously. During the market clearing process, local constraints $\hat{\mathcal{D}}_k$ are not exposed to the operator, which is one advantage of our mechanism.

Thinking of the microgrid operator as a special participant, problem \eqref{eq:sharing-pro} for all $k$ in $\mathcal{I} \cup \mathcal{J}$ and \eqref{eq:sharing-oper} constitute a generalized Nash game \cite{facchinei2010generalized}. The game consists of the following elements:

1) a set of players $\mathcal{K}$, including $I$ consumers indexed by $k=1,...,I$, $J$ prosumers indexed by $k=(I+1),...,(I+J)$ and a microgrid operator indexed by $k=I+J+1$; 

2) action sets $\mathcal{X}_k(\lambda_k):=\{(d_k,b_k)| \; \eqref{eq:sharing-pro.2}-\eqref{eq:sharing-pro.3} \; \mbox{are satisfied}.\}$ for all $k=1,...,(I+J)$, $\mathcal{X}_{I+J+1}(d,b):=\{\lambda | \; \eqref{eq:sharing-oper.2}-\eqref{eq:sharing-oper.3} \; \mbox{are satisfied}.\}$, and action space $\mathcal{X}=\prod_{k} \mathcal{X}_k$; 

3) cost functions $\Gamma_k(\lambda_k):=f_k(d_k)+\lambda_k(-a\lambda_k+b_k)$ for all $k=1,...,(I+J)$, and $\Gamma_{I+J+1}:=\sum_{k \in \mathcal{I} \cup \mathcal{J}} \lambda_k^2$.

The strategies of each user $k \in \mathcal{I} \cup \mathcal{J}$ are its elastic demand $d_k$ and bid $b_k$; the strategy of the microgrid operator is the sharing prices $\lambda_k,\forall k \in \mathcal{I} \cup \mathcal{J}$. When making their decisions, the operator only needs to know the bids $b_k,\forall k \in \mathcal{I} \cup \mathcal{J}$ submitted by users; Each user only needs to know the price $\lambda_k$ given by the operator. Under a market environment, it is reasonable to assume that the microgrid operator can get these bids and the users can get the prices \cite{david2000strategic}. Since the transferred information, i.e. the bids $b_k$ and prices $\lambda_k$, are scalar numbers instead of complex messages, the communication is efficient. For simplicity, we use $\mathcal{G}(w)=\{\mathcal{K},\mathcal{X},\Gamma\}$ to denote the sharing game in an abstract form. Then we define a generalized Nash equilibrium (GNE) as below.
\begin{definition}(Generalized Nash Equilibrium)
A profile $(d^*,b^*,\lambda^*) \in \mathcal{X}$ is a \emph{generalized Nash equilibrium} (GNE) of the energy sharing game $\mathcal{G}(w)=\{\mathcal{K},\mathcal{X},\Gamma\}$ if 
\begin{align}
    (d_k^*,b_k^*) = \mbox{argmin}_{d_k,b_k} \;\{ \Gamma_k(\lambda_k^*), \forall (d_k,b_k) \in \mathcal{X}_k(\lambda_k^*)\}
\end{align}
holds for $k=1,...,(I+J)$, and
\begin{align}
    \lambda^* = \mbox{argmin}_{\lambda} \; \{ \Gamma_{I+J+1}, \forall \lambda \in \mathcal{X}_{I+J+1} (d^*,b^*)\}
\end{align}
\end{definition}

Different from a standard Nash game, in $\mathcal{G}(w)=\{\mathcal{K},\mathcal{X},\Gamma\}$, every player's action set depends on other players' strategies. For example, the action set for user $k=1,...,(I+J)$, the $\mathcal{X}_k$, depends on the operator's strategy $\lambda_k$; the action set for the operator, the $\mathcal{X}_{I+J+1}$, depends on all prosumers' strategies $b$. This complicated coupling makes it difficult to search for and analyze the equilibrium.

\subsection{Properties of the GNE}

The proposed energy sharing market can be used in commercial and industrial microgrids \cite{aram2017microgrid}, campus microgrids \cite{hadjidemetriou2018design}, and community microgrids \cite{long2017peer}. In the following, we discuss its performance. Proposition \ref{prop-1} shows the existence of a market equilibrium with social optimal disutility, while the system's flexibility will be quantified in Section III.

\begin{proposition}
\label{prop-1}
The game $\mathcal{G}(w)=\{\mathcal{K},\mathcal{X},\Gamma\}$ has at least one GNE. Moreover, the triplet $(d^*,b^*,\lambda^*)$ is an GNE if and only if $d^*$ is the unique optimal solution of \eqref{eq:central}, and $\lambda_k^*,\forall k \in \mathcal{I} \cup \mathcal{J}$ equal the corresponding dual variables, with $b_k^*=d_k^f+d_k^*+a\lambda_k^*$ for all $k$ in $\mathcal{I}$, $b_k^*=d_k^f+d_k^*-w_k+a\lambda_k^*$ for all $k$ in $\mathcal{J}$.
\end{proposition}

The proof of Proposition \ref{prop-1} can be found in Appendix \ref{apen-1}. It shows the equivalence of the GNE of $\mathcal{G}$ and the optimal solution of problem \eqref{eq:central}, indicating that the proposed energy sharing mechanism is economically efficient: the GNE $(d^*,b^*,\lambda^*)$ for accommodating $w$ has the lowest social total disutility $\sum_{k \in \mathcal{I} \cup \mathcal{J}} f_k(d_k^*)$.

\textbf{Remark on relevant works:} While both this paper and \cite{dvorkin2019electricity} provide a centralized counterpart for computing the market equilibrium, the problems considered and model configurations are different. Reference \cite{dvorkin2019electricity} investigated a case where different agents (price-setting agent, producer, and consumer) have different information on the probability distribution of renewable uncertainty. Even in the equilibrium problem (2) of \cite{dvorkin2019electricity}, there is still such kind of information asymmetry. In this paper, we consider a case where different agents have access to different constraint sets. To ensure information privacy, an energy sharing mechanism is proposed, resulting in an equilibrium problem where each agent (microgrid operator or user) makes decisions according to its private information, which is known exactly to it. Reference \cite{margellos2017distributed} proposes a proximal minimization based algorithm to deal with the distributed constrained optimization. The constraint $X_i$ is imposed on each local participant $i$; each agent receives information only from its out-neighbours. In this paper, we take into account the coupling constraints $\tilde{\mathcal{P}}$ which is imposed on all users. Besides, instead of exchanging information with neighbors, user only communicates with the microgrid operator.

\textbf{Remark on possible extensions:} (i) Although the above analyses are based on the example of massive users in a stand-alone microgrid, the proposed model and mechanism can be extended to various cases with the main properties remained. For example, if we allow $d_k$ to be negative, the model can incorporate the case with electricity market. To be specific, when $d_k>0$ at the optimal point (or GNE), it means user $k \in \mathcal{I} \cup \mathcal{J}$ adjusts its elastic demand to $d_k$ or when $d_k$ exceeds the upper bound of elastic demand, it also sells electricity to the power market; when $d_k<0$ at the optimal point (or GNE), the user $k$ adjusts its elastic demand to zero, and besides buys $-d_k$ from the power market. (ii) Although there might be some complex units with nonconvex constraints in the microgrid, our approach can still be applied via some convex relaxation techniques, such as those for AC power flow models \cite{farivar2013branch,gan2014convex} and for storage-like devices \cite{li2015storage}.

\subsection{Modified best-response based algorithm}
We prove that the GNE under the proposed mechanism has an appealing property. How to reach such an equilibrium is another important issue. In this paper, a modified best-response (BR) based algorithm (Algorithm 1) is developed. We choose the BR based algorithm as it is one of the most fundamental method in game theory \cite{nash1950equilibrium}. This approach iteratively solves user's problem \eqref{eq:sharing-pro} or the modified market clearing problem \eqref{eq:sharing-oper} with a modified objective function \eqref{eq:modified-function} given other players' strategies until convergence. 

\begin{algorithm}[t]
	\caption{Modified Best-Response}
	\LinesNumbered 
	\KwIn{parameters $a$, $d_k^f,\forall k \in \mathcal{I} \cup \mathcal{J}$, $w_k,\forall k \in \mathcal{J}$, disutility functions $f_k(.),\forall k \in \mathcal{I} \cup \mathcal{J}$}
	\KwOut{generalized Nash equilibrium $(d^*,b^*,\lambda^*)$}
	Initialization: $n=0$, $b^0=0$\; 
	\Repeat{$|b^{n}-b^{n-1}| \le \epsilon$}{
		\emph{Operator}: 
		
		given all bids $b_k^{n},\forall k \in \mathcal{I} \cup \mathcal{J}$, solve problem \eqref{eq:sharing-oper} with a modified objective function:
		\begin{align}
		\label{eq:modified-function}
		\mathop{\min}_{\lambda_k,\forall k \in \mathcal{I} \cup \mathcal{J}} \sum \limits_{\forall k \in \mathcal{I} \cup \mathcal{J}} \lambda_k^2+ \sum \limits_{\forall k \in \mathcal{I} \cup \mathcal{J}} (\lambda_k-\lambda_k^{n})^2
		\end{align}
		to update the price  $\lambda_k^{n+1},\forall k \in \mathcal{I} \cup \mathcal{J}$.
		
		\For{User $k \in \mathcal{I} \cup \mathcal{J}$}{
		 given $\lambda_k^{n+1}$, solves problem \eqref{eq:sharing-pro} to get $d_k^{n+1}$ and $b_k^{n+1}$.}

\emph{Iteration}: $n++$
}
\end{algorithm}

The modified objective function \eqref{eq:modified-function} in Algorithm 1 can smooth the fluctuation of market prices during the bidding process. The value of the market sensitivity $a$ will influence the performance of the modified BR based algorithm. We prove in Appendix \ref{apen-2} that when Condition C1 holds, the bidding process converges.

\noindent \textbf{C1.} The Hessian Matrix of $\sum \nolimits_{k \in \mathcal{I} \cup \mathcal{J}} f_k(d_k)-\frac{1}{2a} \sum \nolimits_{k \in \mathcal{I} \cup \mathcal{J}} (d_k+D_k)^2$  is definite, where $D_k=d_k^f,\forall k \in \mathcal{I}$, $D_k=d_k^f-w_k,\forall k \in \mathcal{J}$.

\begin{figure}[t]
	\centering
	\includegraphics[width=0.8\columnwidth]{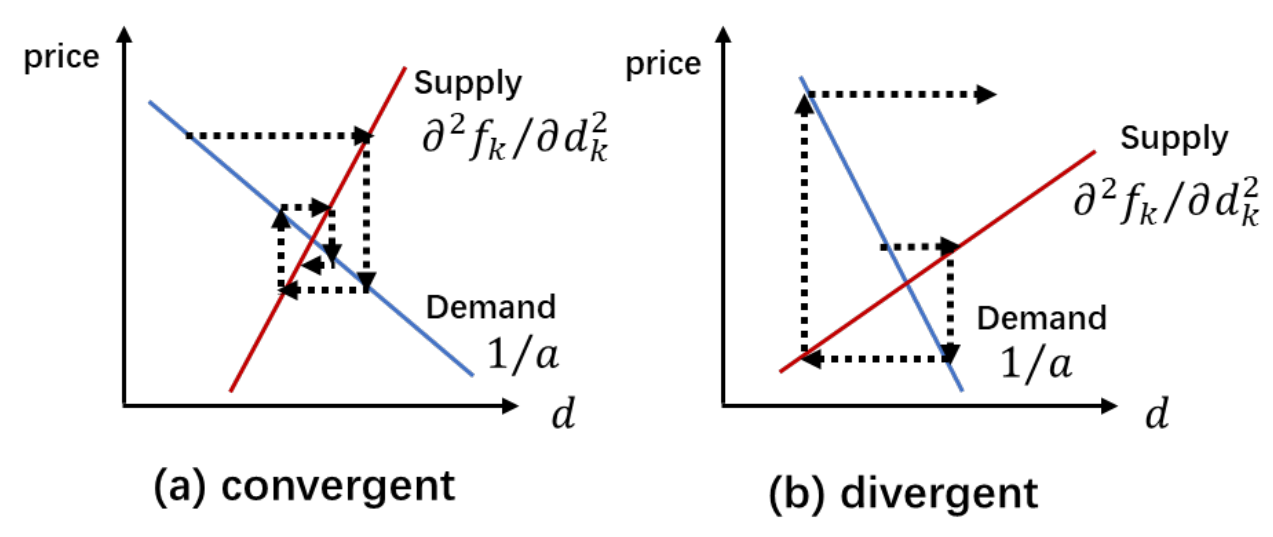}
	\caption{Economic intuition behind Condition C1.}
	\label{fig:intuition}
\end{figure}

Condition C1 can be equivalently written as $1/a < \mathop{\min}\{\partial^2 f_k/\partial d_k^2,\forall k \in \mathcal{I} \cup \mathcal{J}\}$. Specially, if the disutility function is quadratic, i.e. $f_k(d_k)=\alpha_k^1(d_k)^2+\alpha_k^2(d_k)$ with given constants $\alpha_k^1,\alpha_k^2>0$, Condition C1 holds if and only if $a>\mathop{\max}\{\frac{1}{2\alpha_k^1},\forall k \in \mathcal{I} \cup \mathcal{J}\}$. We try to give an economic interpretation of Condition C1 in Fig. \ref{fig:intuition}. $\partial f_k /\partial d_k$ is the marginal disutility of user $k \in \mathcal{I} \cup \mathcal{J}$ and can be regarded as the supply curve of load adjustment, where $\partial^2 f_k/\partial d_k^2$ is the slope of this supply curve. The demand for load adjustment comes from the sharing market, and since $q_k=-a\lambda_k+b_k$, the slope of demand curve is $-1/a$. From Fig. \ref{fig:intuition}, we can find that the bidding process converges if and only if $\partial^2 f_k/\partial d_k^2>1/a,\forall k \in \mathcal{I} \cup \mathcal{J}$, which is also known as the cobweb model in economics. Furthermore, When $a$ is small ($1/a$ is large), the sharing price responses more quickly to the changing bids, and therefore, the algorithm can reach the equilibrium faster.

\section{Flexibility under Energy Sharing}

In this section, we show that energy sharing has the same flexibility as centralized dispatch. The ``absorbable region'' is proposed to characterize the flexibility.

\subsection{Dispatchable region of centralized operation}

In the centralized operation problem (\ref{eq:central}), the renewable outputs $w_j,\forall j \in \mathcal{J}$ are volatile and uncertain. In real-time, the microgrid operator adjusts the $d_k,\forall k \in \mathcal{I} \cup \mathcal{J}$ to ensure supply-demand balance. One critical issue is the system's ability to accommodate uncertainty, which can be quantified by the ``\emph{dispatchable region}'' proposed in \cite{wei2014dispatchable}.
\begin{figure}[h]
        \centering
        \includegraphics[width=0.7\columnwidth]{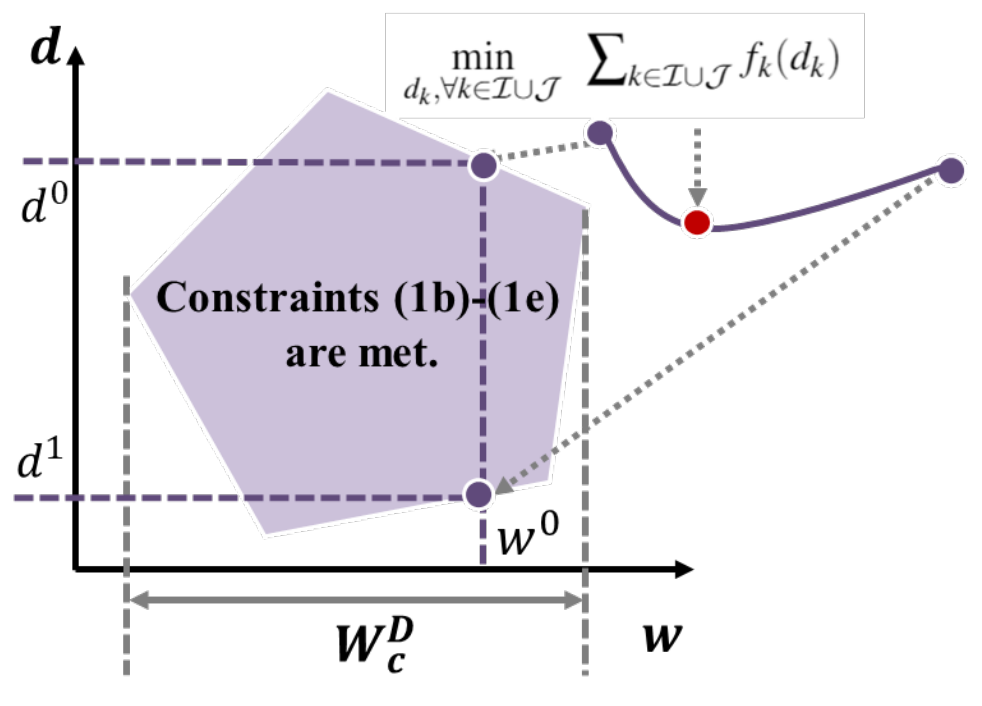}
        \caption{Illustrative diagram of the dispatchable region.}
        \label{fig:DR}
\end{figure}
\begin{definition}(Dispatchable Region \cite{wei2014dispatchable})
The \emph{dispatchable region} under centralized operation is a set of renewable outputs such that at least one feasible dispatch solution exists:
$$
W_c^{D}=\{w| \; \exists ~ d : \eqref{eq:central.2}-\eqref{eq:central.5} \; \mbox{are satisfied}.\}
$$
\end{definition}
Depicting the flexibility of a power grid under centralized operation with the dispatchable region has been widely studied \cite{wei2015real,wan2016maximum}. An illustrative diagram is given in Fig. \ref{fig:DR}.
The horizontal axis of the figure is the uncertain factor $w$, and the vertical axis is the elastic demand $d$. The light purple region is the feasible region of problem \eqref{eq:central} characterized by constraints \eqref{eq:central.2}-\eqref{eq:central.5}. Projecting it onto the horizontal axis we can obtain a range of $w$, the ``dispatchable region'' $W_c^D$. For any given point $w^0$ within $W_c^D$, we can always find a corresponding feasible range for $d$, i.e. $[d^0,d^1]$. The operator minimizes the objective function $\sum \nolimits_{k \in \mathcal{I} \cup \mathcal{J}} f_k(d_k)$ over $[d^0,d^1]$ at the red point.

In recent years, innovative approaches emerge for running the power system in a smarter and more scalable way, which calls for methods to quantify the flexibility under those approaches. In this paper, we generalize the conventional dispatchable region to the \emph{``absorbable region"} under any certain schemes for characterizing the flexibility under our energy sharing mechanism.

\subsection{Absorbable Region Under Energy Sharing}
The dispatchable region is a useful tool in describing the system's flexibility but limited to centralized operation. First, we raise the concept of ``absorbable region'' which generalizes the ``dispatchable region'' for flexibility evaluation.

\begin{definition}(Absorbable Region) The \emph{absorbable region} of a microgrid under \emph{a certain scheme} is a set of renewable outputs such that at least one feasible \emph{strategy} exists under that scheme.
\end{definition}

The ``absorbable region'' differs and improves from the ``dispatchable region'' in two ways: Firstly, it is designed for any given scheme instead of just the centralized operation. Secondly, the renewable output is ``absorbable'' if and only if there exists a strategy, which could not only be a dispatch order but also a market equilibrium: When the ``certain scheme'' refers to centralized operation, the absorbable region degenerates to the dispatchable region. If under ``certain scheme'' every user makes decision individually, then the renewable output is absorbable when there is a feasible self-sufficient strategy, i.e. $(w_k-d_k^f) \in \hat D_k,\forall k \in \mathcal{J}$. When the ``certain scheme'' refers to a market mechanism, ``absorbable'' is equivalent to the existence of a market equilibrium. Specially, the absorbable region under the proposed energy sharing mechanism is given by
\begin{align}
W_s^D=\{w\;| \; \mbox{an GNE} \; (d^*,b^*,\lambda^*) \; \mbox{for the game} \; \mathcal{G}(w) \; \mbox{exists}.\} \nonumber
\end{align}

The proposed absorbable region can accommodate various scenarios, however, it may also encounter complicated situations, such as conflicting interests among stakeholders and time or spatial coupling constraints. Though fruitful works have been conducted for computing the dispatchable region under centralized operation, the proposed methods cannot be applied directly to the characterization of the absorbable region in general cases due to the above complexity. In the following, we show that our proposed energy sharing mechanism is flexibility-retained, which facilities the computation of its absorbable region:

According to Proposition \ref{prop-1}, given a $w$, if problem \eqref{eq:central} is feasible, meaning that $w \in W_c^D$, then we can always construct a GNE for the game $\mathcal{G}(w)$, so that $w \in W_s^D$. Conversely, if $w \in W_s^D$, i.e. there exists a GNE for the game $\mathcal{G}(w)$, then the $d^*$ is exactly the optimal solution (of course also a feasible one) of the problem \eqref{eq:central}, indicating $w \in W_c^D$. Thus, we have $W_c^D=W_s^D$. The proposed energy sharing mechanism can achieve the same flexibility as the centralized dispatch. We can use the simpler equivalent model \eqref{eq:central} for computing the absorbable region of the proposed energy sharing mechanism.

\subsection{Linear programming based projection algorithm}

Here, we use a common power system model with capacity constraints ($\hat{\mathcal{D}}_k$ is chosen as $[\underline{d}_k,\overline{d}_k]$) and network constraints ($\tilde{\mathcal{P}}$ is chosen as the DC power flow limits) as an example. In such a case, constraints \eqref{eq:central.2}-\eqref{eq:central.5} are
\bsq
\label{eq:dispatchable-region}
\begin{align}
        & \sum \limits_{k \in \mathcal{J}} w_k = \sum \limits_{k \in \mathcal{I}\cup \mathcal{J}} d_k + \sum \limits_{k \in \mathcal{I} \cup \mathcal{J}} d_k^f \\
        & \underline{d}_k \le d_k \le \overline{d}_k,~\forall k \in \mathcal{I} \cup \mathcal{J} \\
        & -F_l \le \sum \limits_{k \in \mathcal{J}} \pi_{kl}^w w_k - \sum \limits_{k \in \mathcal{I} \cup \mathcal{J}} \pi_{kl}^d (d_k^f+d_k) \le F_l,\forall l \in \mathcal{L}
\end{align}
\esq
where $\mathcal{L}$ is the set of lines, $F_l$ is the power flow limit for line $l \in \mathcal{L}$, and $\pi_{kl}^w, \pi_{kl}^d$ are the line flow distribution factors. Since all constraints in \eqref{eq:dispatchable-region} are linear, it can be represented in a compact form:
\begin{align}
        \psi(w)= \{d~|~ Ad+Bw \le c\}
\end{align}
Here, $d$ is a collection of $d_k,\forall k \in \mathcal{I} \cup \mathcal{J}$, and $w$ is a collection of $w_k,\forall k \in \mathcal{J}$. According to Proposition \ref{prop-1}, the absorbable region is defined as 
\begin{align}
        W_s^D = W_c^D= \{w ~|~ \psi(w) \ne \emptyset\}
\end{align}

For a given $w$, we can check whether $\psi(w)$ is empty by the following problem:
\bsq
\label{eq:feasibility-check}
\begin{align}
        \mathop{\min}_{d,z} ~ & 1^T z \\
        \mbox{s.t.}~ & Ad-Iz \le c-Bw : u \\
        ~ &  z \ge 0
\end{align}
\esq
where $z$ is a slack variable and $u$ is the dual variable. It is easy to see that $\psi(w) \ne \emptyset$ if and only if the optimal value of \eqref{eq:feasibility-check} is zero. The dual problem of the linear program \eqref{eq:feasibility-check} is 
\bsq
\begin{align}
        \mathop{\max}_{u}~ & u^T(c-Bw) \\
        \mbox{s.t.}~ & A^Tu=0, ~ -1\le u \le 0
\end{align}
\esq
Denote $\mathcal{U}:=\{u~|~ A^Tu=0,-1 \le u \le 0\}$, then $\psi(w) \ne \emptyset$ is equivalently to
\begin{align}
u^T(c-Bw) \le 0,\forall u \in \mathcal{U}
\label{eq:inf-enumer}
\end{align} 
The set $\mathcal{U}$ is a closed convex set, so the above condition can be further simplified to
\begin{align}
\label{eq:feasibility-eq}
        u^T(c-Bw) \le 0,\forall u \in \mbox{vert}(\mathcal{U})
\end{align}

Consider $w$ as the variable, polyhedron \eqref{eq:feasibility-eq} is the condition which $w$ must meet in order to ensure the non-emptiness of $\psi(w)$, i.e.:
\begin{equation}
\label{eq:closed-form}
W_s^D= \{w ~|~ u^T(c-Bw) \le 0,\forall u \in \mbox{vert}(\mathcal{U})\}
\end{equation}
Though $W_s^D$ can be represented as an explicit polyhedron, it is still difficult to locate all the vertices of $\mathcal{U}$ with a high-dimension. The projection algorithm (Algorithm 2) is developed to generate the absorbable region $W_s^D$ efficiently, in order to exhibit that the distributed method possesses the same flexibility as the centralized method. This projection algorithm is not actually executed in the market clearing procedure and consumes real-time computational resources.

\begin{algorithm}[t]
        \caption{Linear Programming Based Projection}
        \LinesNumbered 
        \KwIn{initial set $W_{temp}=\{w|~w \ge 0\}$}
        \KwOut{output dispatchable region $W_s^D$}
        Update vert($W_{temp}$)\; 
        \For{$w \in \mbox{vert}$($W_{temp})$}{
                solve problem
                $$\mathop{\max}_{u}~ u^T(c-Bw) ~~ \mbox{s.t.} ~u \in \mathcal{U}   $$
        and denote the optimal solution as $u^*$, the optimal value as $r^*$. Let $r_{max}=\mathop{\max}\{r_{max},r^*\}$ and update $u_{max}$ as the corresponding $u^*$.
        }
        \eIf{$r_{max}=0$}{
                        let $W_s^D=W_{temp}$\;
                }{
                        add a cutting plane $(u_{max})^TBw \ge (u_{max})^Tc$ in $W_{temp}$\;
                        go to 1\;
                }
\end{algorithm}

The intuition behind Algorithm 2 is as follows. We know that $w \in W_s^D$ if and only if \eqref{eq:feasibility-eq} holds; $w \notin W_s^D$ indicates that $(u_{max})^T(c-Bw)>0$ for some $u_{max} \in \mbox{vert}(\mathcal{U})$. Therefore, we first  initiate a large enough polyhedron $W_{temp}$ which contains $W^D_s$, and then test if (\ref{eq:inf-enumer}) is met by solving the linear program in step 3. If not, the hyperplane 
\begin{align}
        (u_{max})^Tc= (u_{max})^TBw 
\end{align}
gives the boundary of $W_s^D$, which is added to $W_{temp}$. The algorithm converges  once  (\ref{eq:inf-enumer}) is met. Algorithm 2 is more efficient than the mixed-integer programming based algorithm in \cite{wei2014dispatchable} and does not need a big-M parameter.

\textbf{Remark on convergence of Algorithm 2}: The convergence criteria of Algorithm 2 is that \eqref{eq:inf-enumer} is met for all $w$ in $W_s^D$. The region $W_s^D$ is a polyhedron with a finite number of plane faces and the closed-form in \eqref{eq:closed-form} based on enumerating the vertices of $\mathcal{U}$. However, when the dimension of $\mathcal{U}$ is high, it is challenging to list all its extreme points. Actually, most of these extreme points only exert redundant constraints in $W_s^D$ which can be removed. Algorithm 2 avoids these redundant constraints by adaptively selecting the vertices in $\mathcal{U}$ and adds a hyperplane to give the boundary of $W_s^D$  in each iteration. Since there are a finite number of vertices of $\mathcal{U}$, the algorithm can always converge. When there are many constraints in the problem \eqref{eq:dispatchable-region}, the polyhedron $W_s^D$ is likely to have many facets, and thus needs more iterations to characterize. Nonetheless, as the procedures in the Algorithm 2 only require solving linear programs, the algorithm is efficient.

\section{Case Studies}
Numerical examples are tested to validate the theoretical results and show the effectiveness of the proposed algorithm.

\subsection{Simple case with two groups of prosumers}
\begin{figure}[h]
        \centering
        \includegraphics[width=0.7\columnwidth]{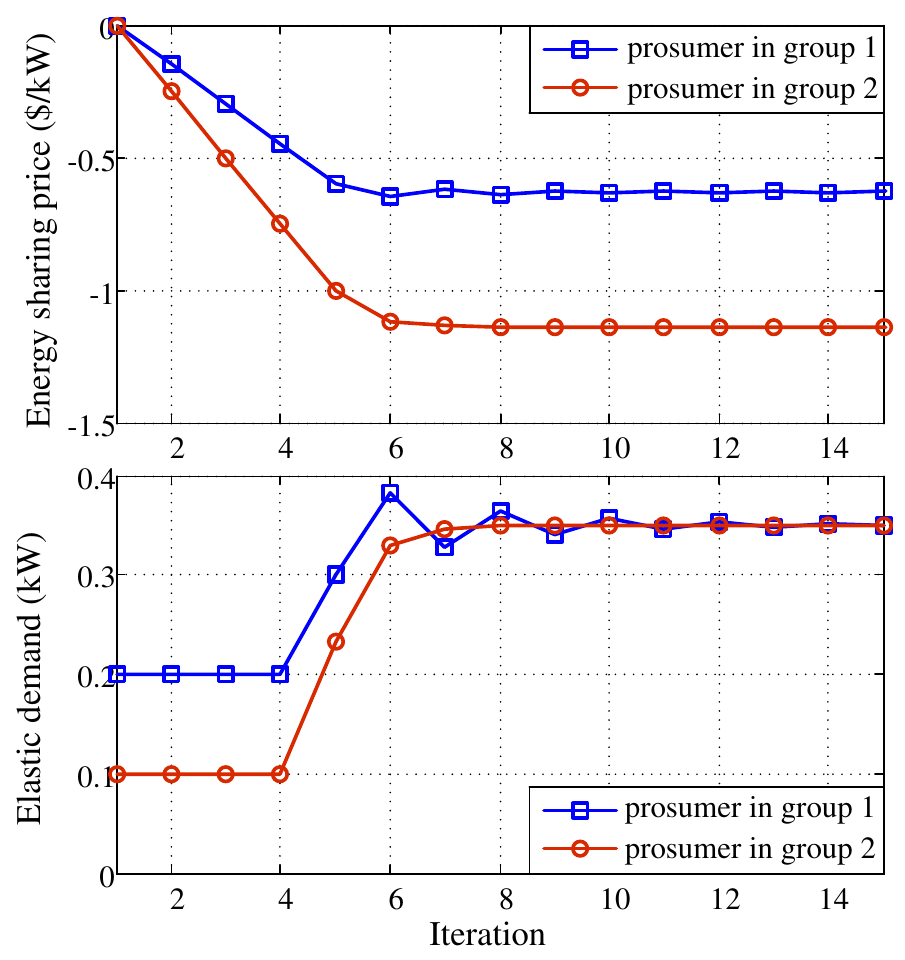}
        \caption{Sharing prices and optimal strategies in each iteration.}
        \label{fig:convergence}
\end{figure}

First, we verify Proposition 1 by a simple example with two groups of prosumers. The first group consists of prosumers $k \in \mathcal{J}_1=\{1,...,100\}$, and the second group consists of prosumers $k \in \mathcal{J}_2=\{101,...,200\}$. Each group is made up of identical prosumers. We adopt quadratic disutility functions $f_k(d_k):=\alpha_k^1 (d_k)^2+ \alpha_k^2 d_k,\forall k \in \mathcal{J}_1 \cup \mathcal{J}_2$ with the parameters given in Table \ref{tab:parameter}. Two groups of prosumers are connected by a line whose flow limit is $10 \; \mbox{kW}$. The sensitivity $a$ is set to 1 kW/\$. The modified best response based algorithm is used for reaching the GNE, and the changes of energy sharing prices and elastic demands are plotted in Fig. \ref{fig:convergence}.
\begin{table}[h]
        \renewcommand{\arraystretch}{1.3}
        \renewcommand{\tabcolsep}{1em}
        \centering
        \caption{Parameters of prosumers}
        \label{tab:parameter}
        \begin{tabular}{cccccc}
                \hline 
                Prosumer & $\alpha_k^1$ & $\alpha_k^2$ & $w_k$ & $d_k^f$ & $d_k$\\
                 group & ($\$/\mbox{kW}^2$) & ($\$/\mbox{kW}$) & ($\mbox{kW}$) & ($\mbox{kW}$) & ($\mbox{kW}$)\\
               \hline
                1 & 0.30 & 0.42 & 1.25 & 1.00 & [0.2,0.5]\\
               2 & 0.60 & 0.72 & 1.75 & 1.30 & [0.1,0.6] \\
                \hline
        \end{tabular}
\end{table}

From Fig. \ref{fig:convergence} we can find that both the energy prices and prosumers' strategies converge. At GNE, we have $d_1=0.35 \mbox{kW}$ and $d_2=0.35 \mbox{kW}$, which equals the optimal solution of problem \eqref{eq:central} with the same parameters. So the proposed energy sharing mechanism achieves social optimum.

\begin{figure}[h]
	\centering
	\includegraphics[width=1.0\columnwidth]{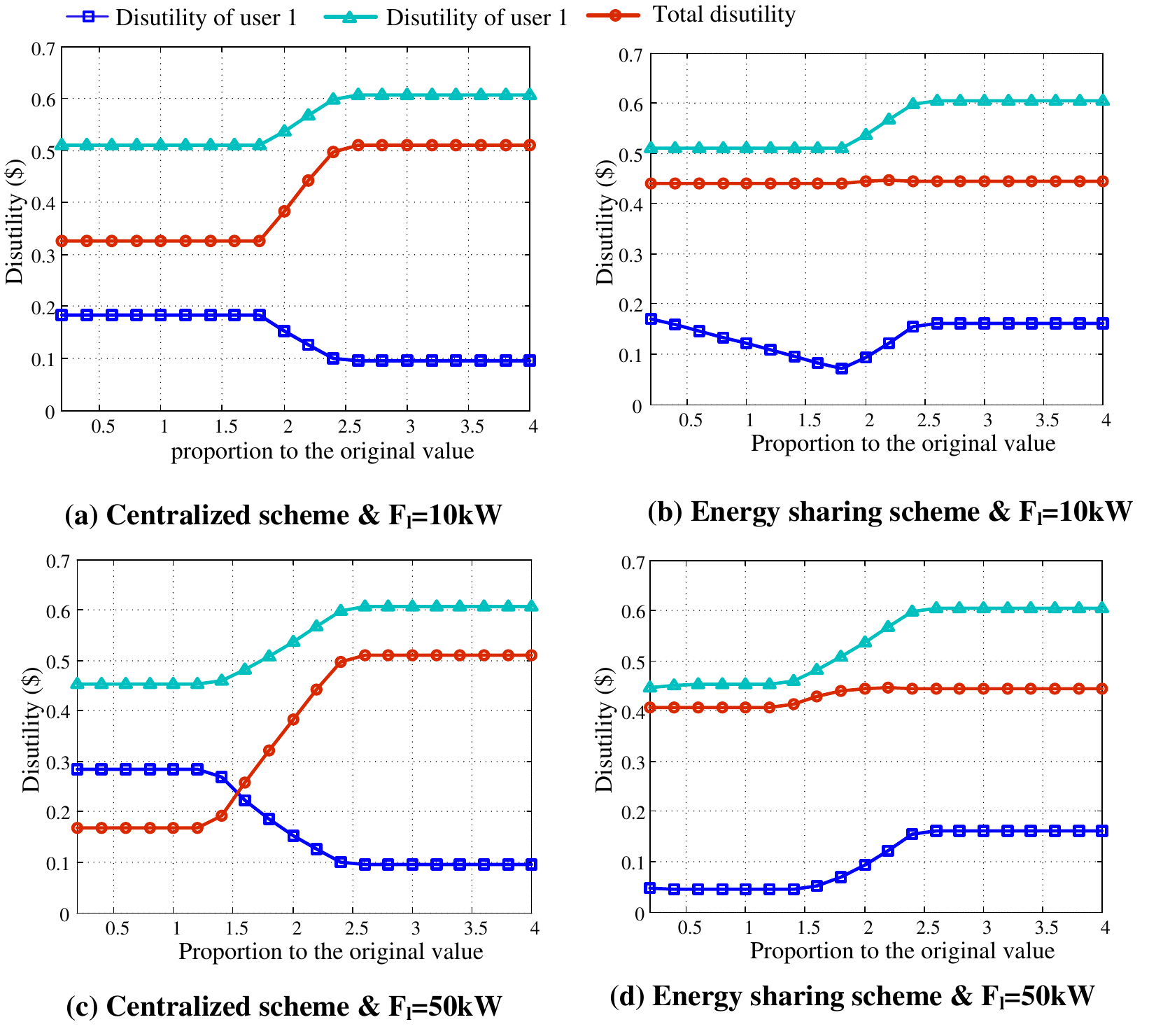}
	\caption{Impact of information asymmetry.}
	\label{fig:asymmetric-info}
\end{figure}

We further reveal the impact of information asymmetry by testing how the misrepresentation of $\alpha_k^1$ and $\alpha_k^2$ will influence the users’ disutilities and social total disutility. We change $\alpha_1^1$ and $\alpha_1^2$ from 0.2 to 4 times of their original values $\bar \alpha_1^1=0.30 \mbox{\$/kW}^2$, $\bar \alpha_1^2=0.42 \mbox{\$/kW}$, and test the case under $F_l=10 \mbox{kW}$ and $F_l=50 \mbox{kW}$, respectively. Results are shown in Fig. \ref{fig:asymmetric-info}. From (a) and (c), we can find that under centralized scheme, user 1 can lower its disutility by reporting larger cost coefficients so that it would be asked to adjust less by the operator. While user 1 can take the advantage of asymmetric information, both user 2’s disutility and the total disutility increase dramatically. When it comes to the proposed energy sharing market: in (b), though user 1 can still obtain a lowest disutility by reporting $1.6\bar \alpha_1^1$ and $1.6\bar \alpha_1^2$, the total disutility remains unchanged; in (d), user 1 has no incentive to misreport because its lowest disutility is with $\bar \alpha_1^1$ and $\bar \alpha_1^2$. The impact of information asymmetry on the total disutility can be lessened via the proposed energy sharing scheme.

\subsection{Five-bus system for flexibility characterization}
\begin{figure}[h]
        \centering
        \includegraphics[width=1.0\columnwidth]{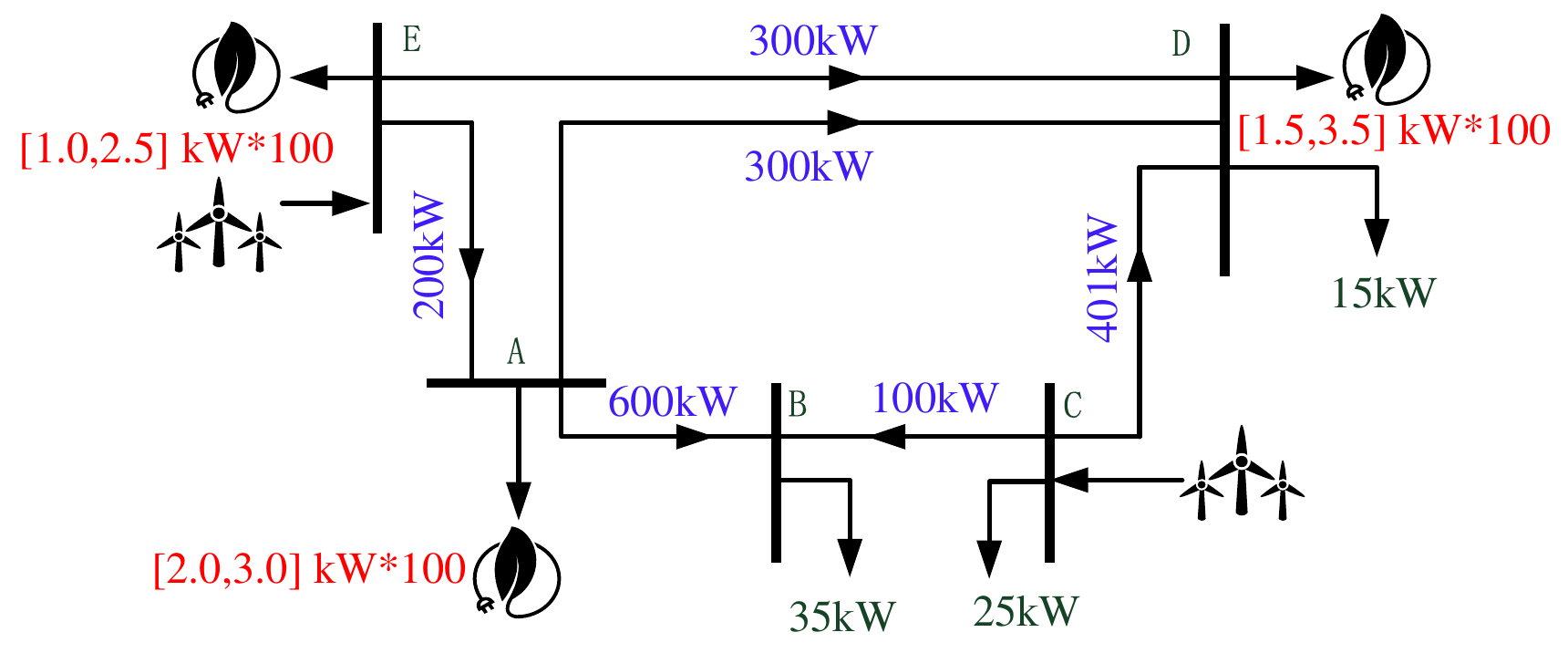}
        \caption{Topology and parameters of the five-bus system.}
        \label{fig:typology}
\end{figure}

\begin{figure}[h]
        \centering
        \includegraphics[width=1.0\columnwidth]{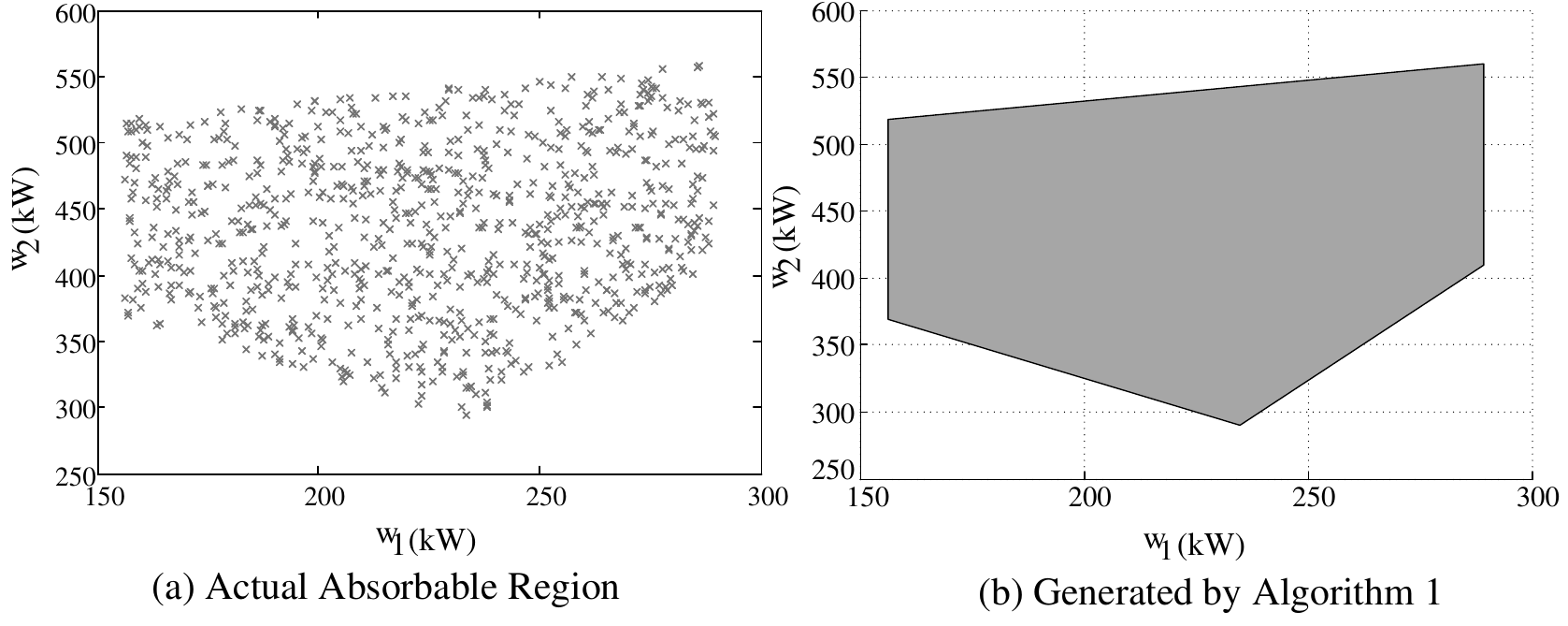}
        \caption{Absorbable region: actual (a) v.s. generated by Algorithm 2 (b).}
        \label{fig:region-PJM}
\end{figure}

Next, the effectiveness of Algorithm 2 for computing the absorbable region is tested via a five-bus system, whose topology and parameters are given in Fig. \ref{fig:typology}. There are three elastic demands, whose range is marked in red; the power flow limits are in blue; the fixed demands are in green. 
We randomly choose renewable output scenarios $(w_1,w_2)$, and mark those with which the problem \eqref{eq:central} is feasible (thus $(w_1,w_2) \in W_s^D$) in Fig.  \ref{fig:region-PJM}(a). The one provided by Algorithm 2 is shown in Fig. \ref{fig:region-PJM}(b). Both regions are exactly the same, demonstrating that Algorithm 2 can successfully identify the boundary of the actual absorbable region. Any renewable output profile $(w_1,w_2)$ inside the grey region is \emph{absorbable}.


\begin{figure}[h]
        \centering
        \includegraphics[width=1.0\columnwidth]{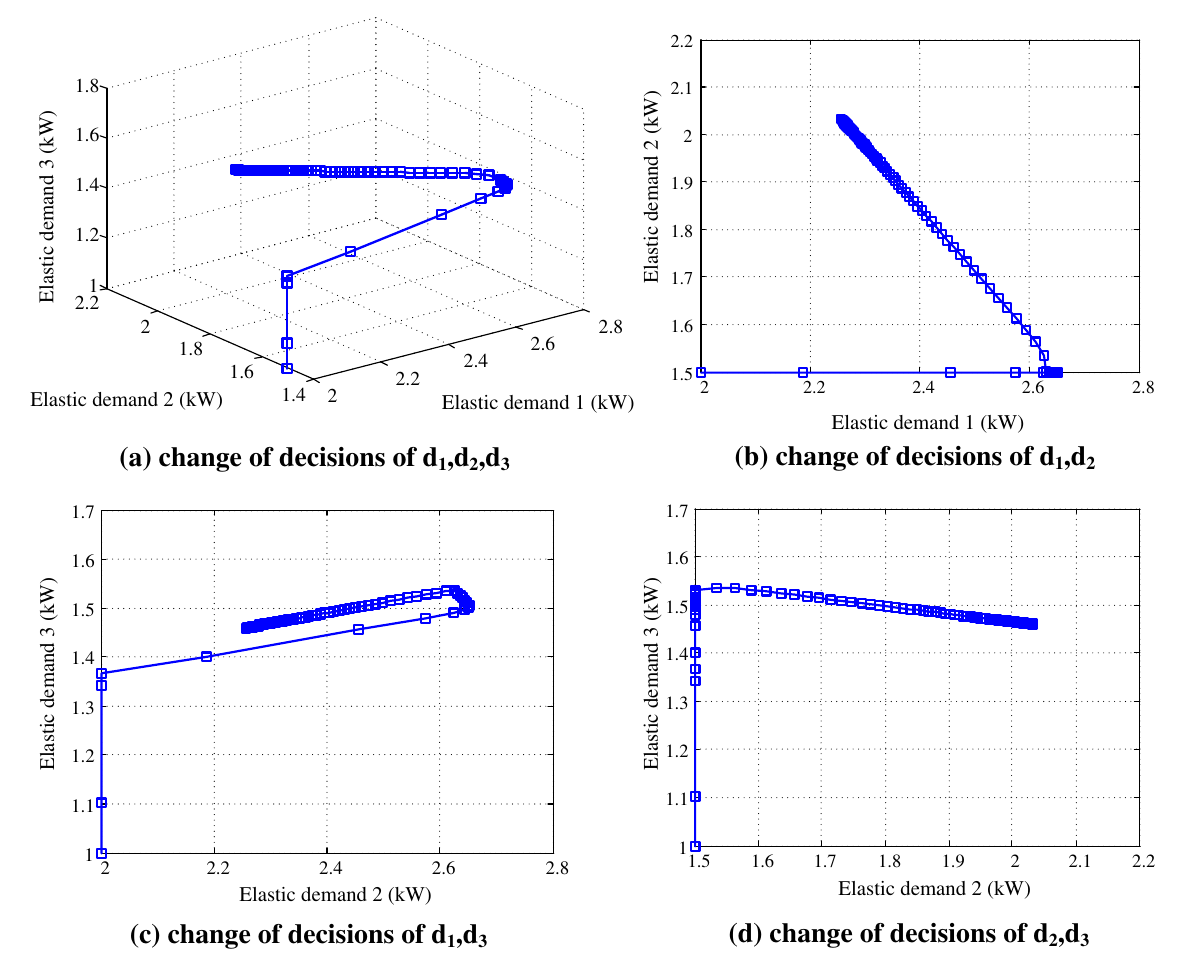}
        \caption{Sequence of elastic demands during iterations in the five-bus system.}
        \label{fig:bid-PJM}
\end{figure}

\begin{figure}[h]
        \centering
        \includegraphics[width=1.0\columnwidth]{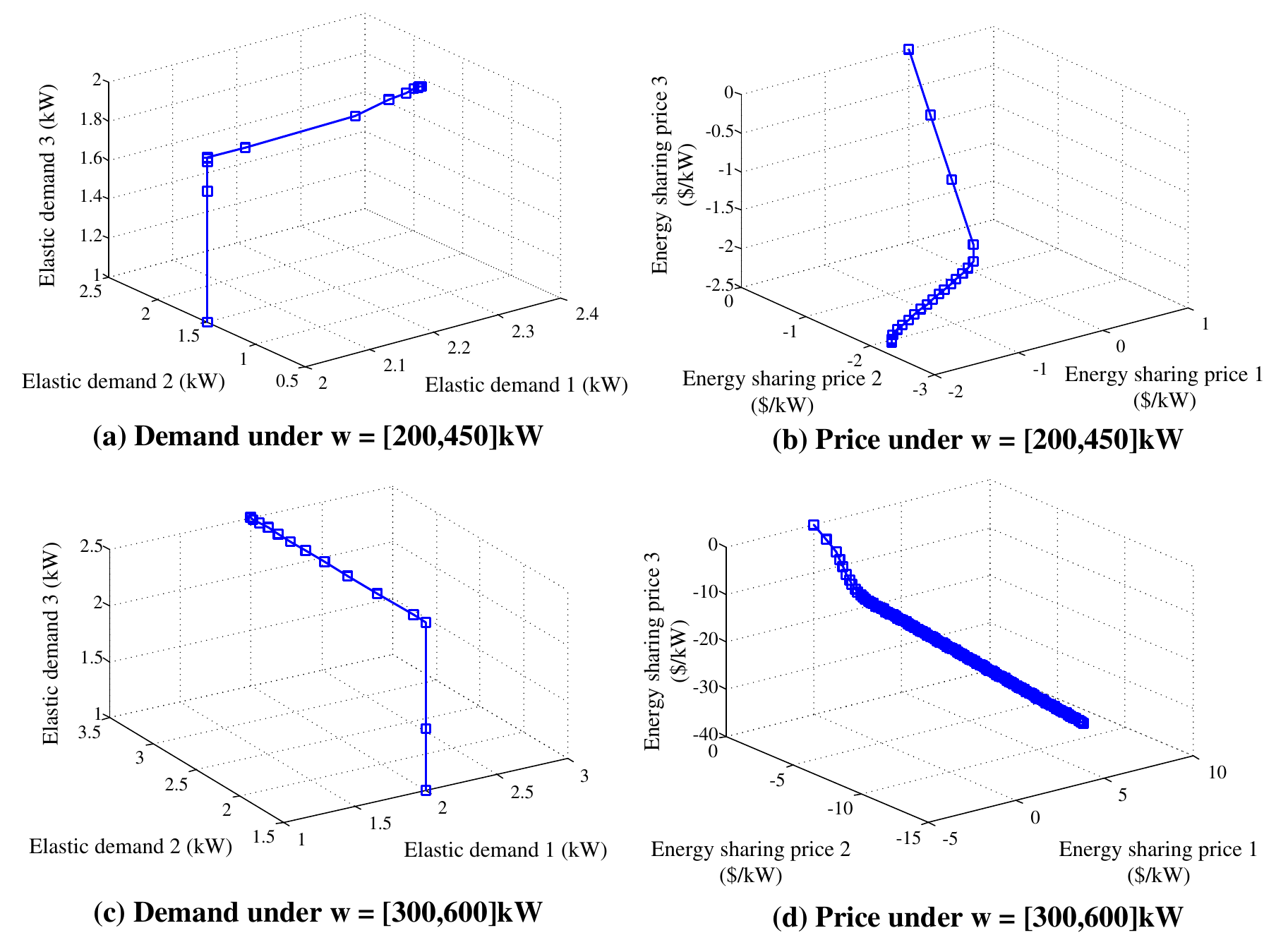}
        \caption{Results under two scenarios for the five-bus system.}
        \label{fig:bid-PJM2}
\end{figure}

Let $w_1=250$ kW, $w_2=400$ kW, and $a=100$ kW/\$, the sequence of the elastic demands generated by Algorithm 1 is shown in Fig. \ref{fig:bid-PJM}(a), and we project it onto $d_1-d_2$,  $d_1-d_3$ and  $d_2-d_3$ planes, respectively, and get Fig. \ref{fig:bid-PJM}(b)-(d). Again, the optimal strategies of three elastic demands converge to $d_1=2.26$ kW, $d_2=2.03$ kW and $d_3=1.46$ kW, which coincide with the optimal dispatch in \eqref{eq:central}. Another two renewable output scenarios are tested with results in Fig. \ref{fig:bid-PJM2}. For the renewable outputs inside the absorbable region, a GNE exists and can be reached by the modified BR based method.

\subsection{Practicability of the proposed mechanism and algorithm}
\begin{figure}[h]
        \centering
        \includegraphics[width=0.8\columnwidth]{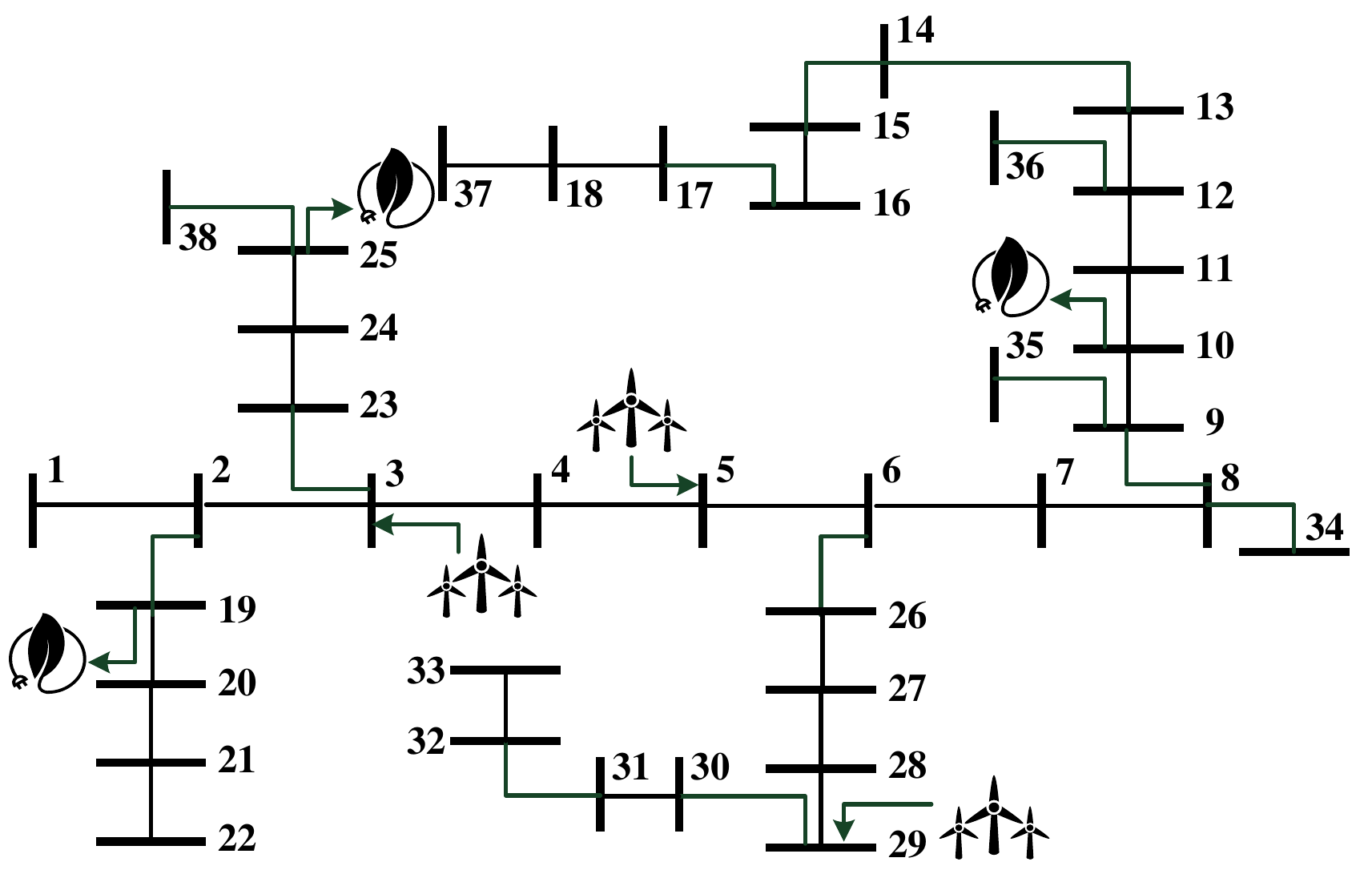}
        \caption{Topology of the 38-bus test system.}
        \label{fig:37bus-topology}
\end{figure}

\begin{figure}[h]
	\centering
	\includegraphics[width=1.0\columnwidth]{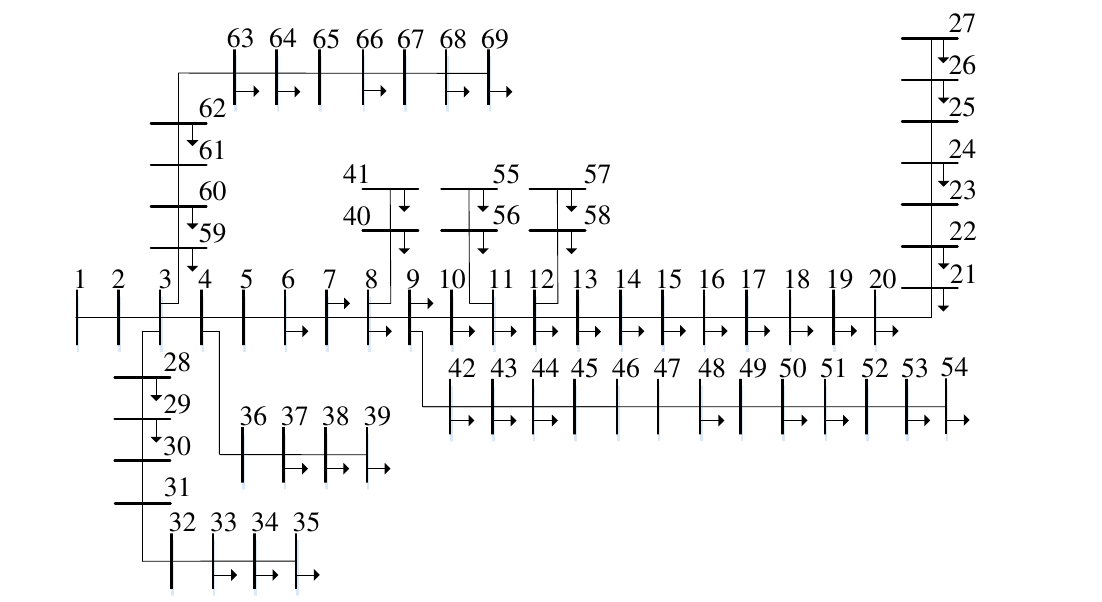}
	\caption{Topology of the 69-bus test system.}
	\label{fig:69bus-topology}
\end{figure}

\begin{figure}[h]
        \centering
        \includegraphics[width=1.0\columnwidth]{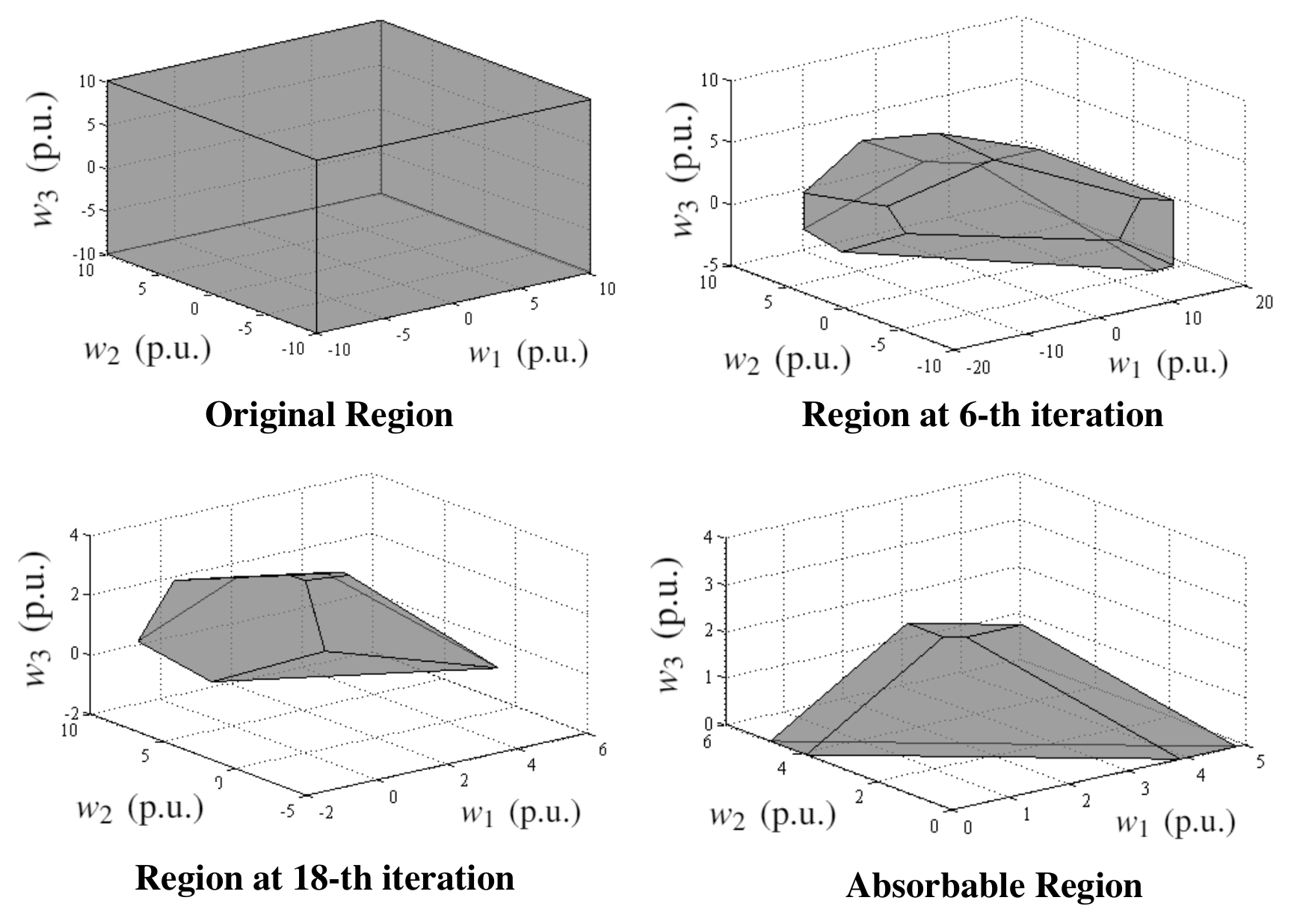}
        \caption{Results when applying Algorithm 2 to the 38-bus system.}
        \label{fig:37bus-AR}
\end{figure}

\begin{figure}[h]
	\centering
	\includegraphics[width=1.0\columnwidth]{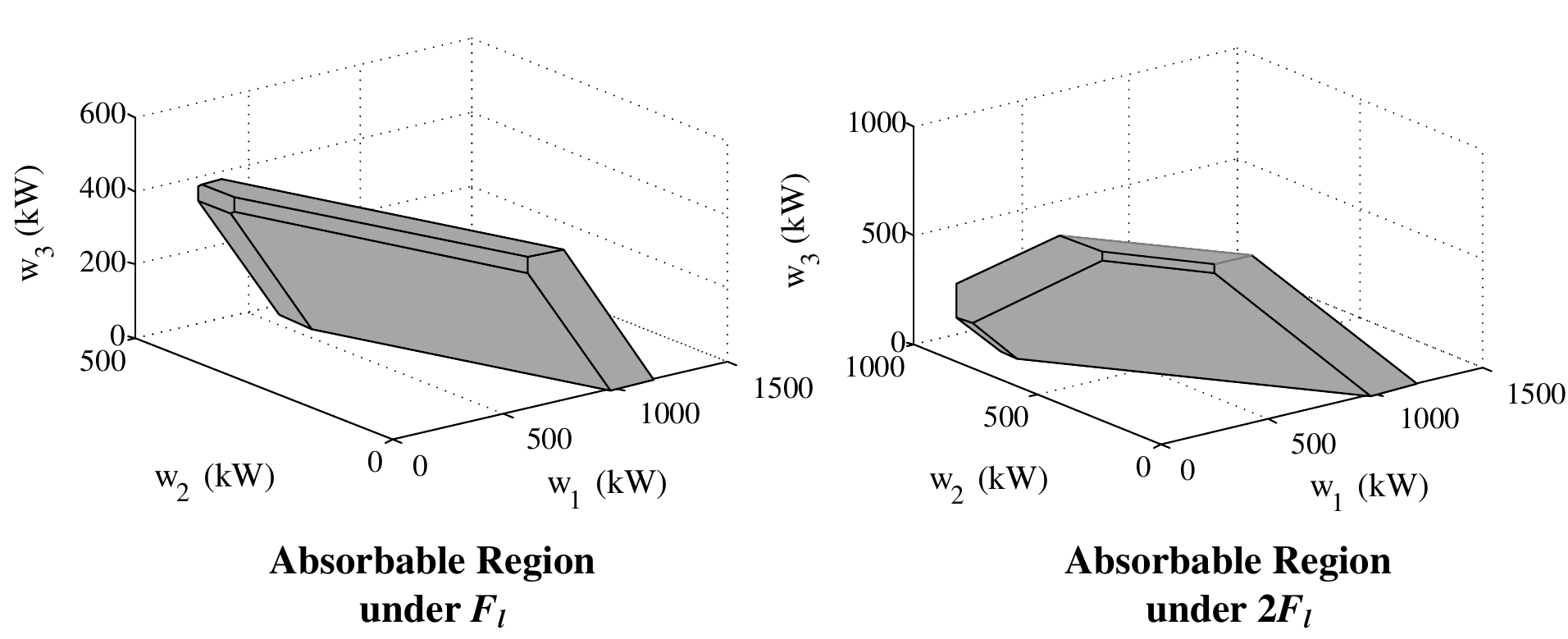}
	\caption{Absorbable region under different power flow limits.}
	\label{fig:69bus-AR}
\end{figure}

A larger case with a modified 38-bus microgrid \cite{singh2009multiobjective} and a modified 69-bus microgrid \cite{kadir2013optimal} are tested to show the scalability of our method. The topology of the test systems are in Fig. \ref{fig:37bus-topology} and Fig. \ref{fig:69bus-topology}, respectively, with other parameters in \cite{Data}. Algorithm 2 is used to generate the absorbable region. The result for 38-bus system is presented in Fig. \ref{fig:37bus-AR} together with some intermediate results. At the beginning, a large enough box that contains the absorbable region is used as the initial polyhedron; then the points outside the absorbable region are gradually removed by  the cutting planes in each iteration. We then generate the absorbable region of the 69-bus system with the flow limits equal to $F_l$ and $2F_l$, respectively, as shown in Fig. \ref{fig:69bus-AR}.
We check that the output absorbable regions are the same as the actual ones, which validates the proposed method. With looser flow limits, the absorbable region enlarges.


\begin{figure}[h]
        \centering
        \includegraphics[width=1.0\columnwidth]{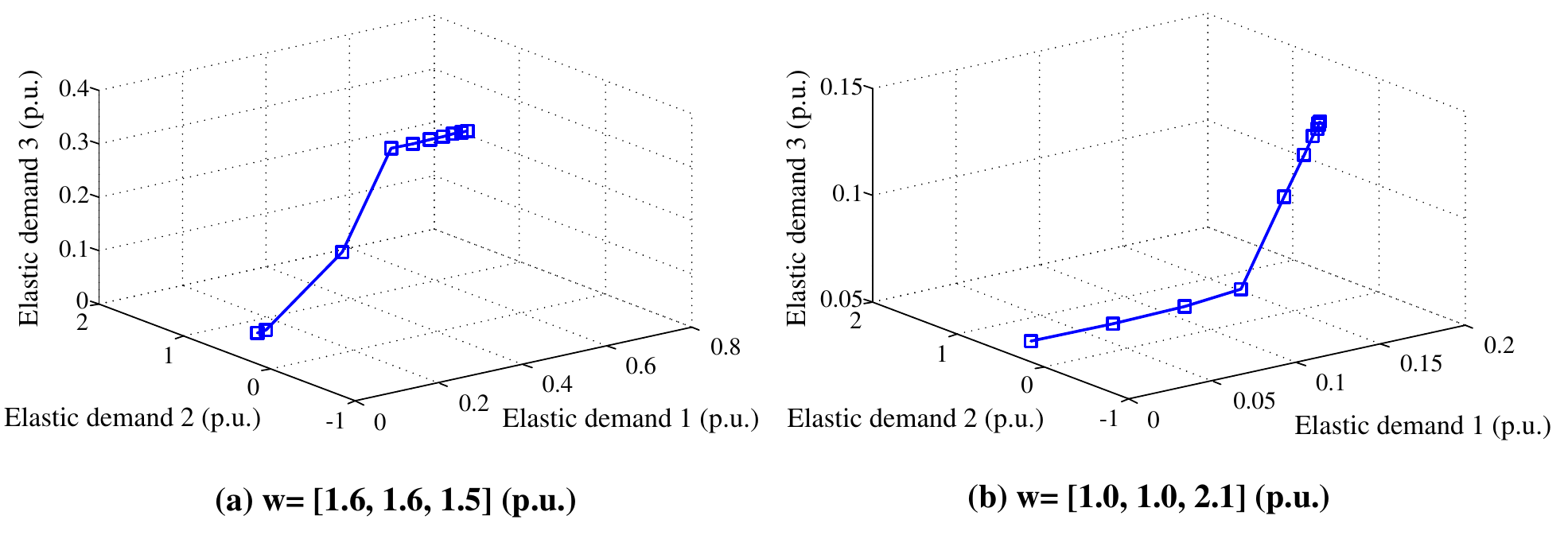}
        \caption{Results under two scenarios for the 38-bus test system.}
        \label{fig:37bus-iteration}
\end{figure}

\begin{figure}[h]
	\centering
	\includegraphics[width=1.0\columnwidth]{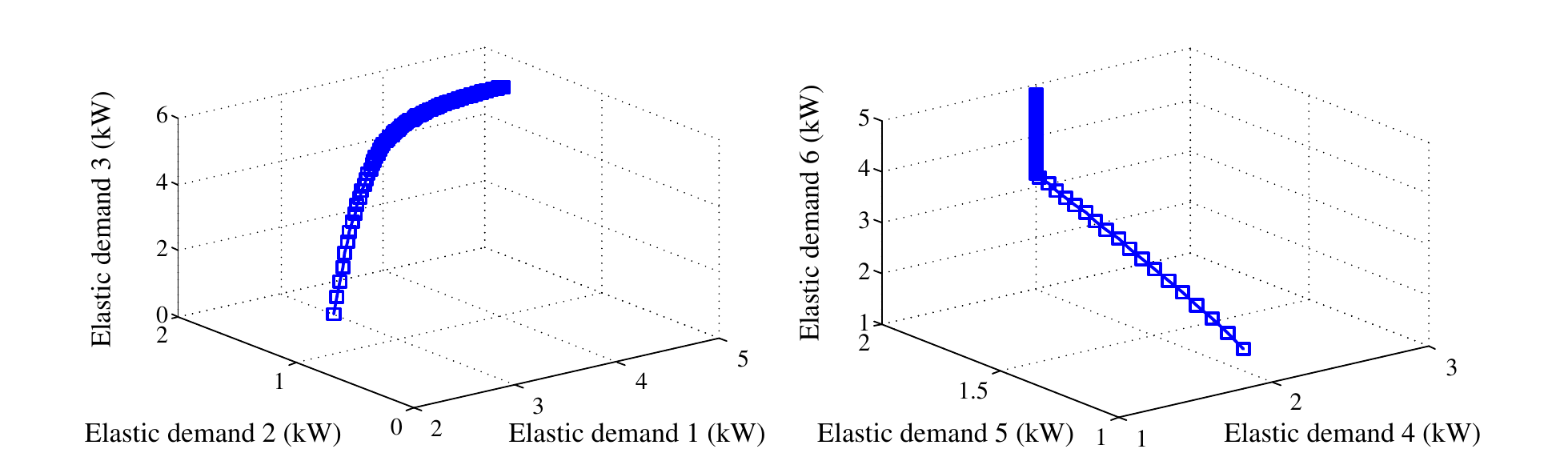}
	\caption{Change of users’ strategies for the 69-bus test system.}
	\label{fig:69bus-iteration}
\end{figure}

We choose two wind output scenarios within the above obtained absorbable region for 38-bus system, and let $w=[1.6, 1.6, 1.5]$ (p.u.) and $w=[1.0, 1.0, 2.1]$ (p.u.), respectively. The change of prosumers' strategies on elastic demands are recorded in Fig. \ref{fig:37bus-iteration}. Under both cases, prosumers' strategies converge within 10 iterations. In addition, the output strategies are $[0.50, 0.20, 0.35]$ (p.u.) and $[0.17, 0.14, 0.13]$ (p.u.), which is the optimal solution of problem \eqref{eq:central}. Similarly, for the 69-bus system, let $w=[300, 450, 400]$ kW, the change of end users' strategies when implementing the modified BR algorithm is given in Fig. \ref{fig:69bus-iteration}, and all of them converge in 40 iterations.


We list the time needed for reaching equilibrium for the simple case in Section IV.A, the five-bus case in Section IV.B, and the 38-bus and 69-bus cases in Section IV.C in TABLE \ref{tab:time}. We can find that all cases converge within 300 seconds, while the scale of the system and the number of consumers/prosumers do not matter much. Our proposed energy sharing mechanism focuses on the real-time market \cite{pei2016optimal}, which will be conducted hourly. Therefore, the time cost of the algorithm is acceptable for its implementation.
\begin{table*}[h]
	\renewcommand{\arraystretch}{1.3}
	\centering
	\caption{Computational time (s) for different cases.}
	\label{tab:time}
	\begin{tabular}{ccccccccc}
		\hline 
		Case & Simple & 5-bus with  & 5-bus with & 5-bus with & 38-bus with  & 38-bus with & 38-bus with  & 38-bus with \\
		  & case & $w=[200,400]$ & $w=[250, 450]$ & $w=[300,600]$ & $w=[1.6,1.6,1.5]$ & $w=[1.0,1.0,2.1]$ & $w=[300,450,400]$ & $w=[400,350,300]$ \\
		Time & 24.85 & 238.32 & 65.85 & 47.55 & 23.61 & 32.20 & 143.14 & 132.83\\
		\hline
	\end{tabular}
\end{table*}

\subsection{Factors influencing the performance of algorithms}
\begin{figure}[h]
	\centering
	\includegraphics[width=1.0\columnwidth]{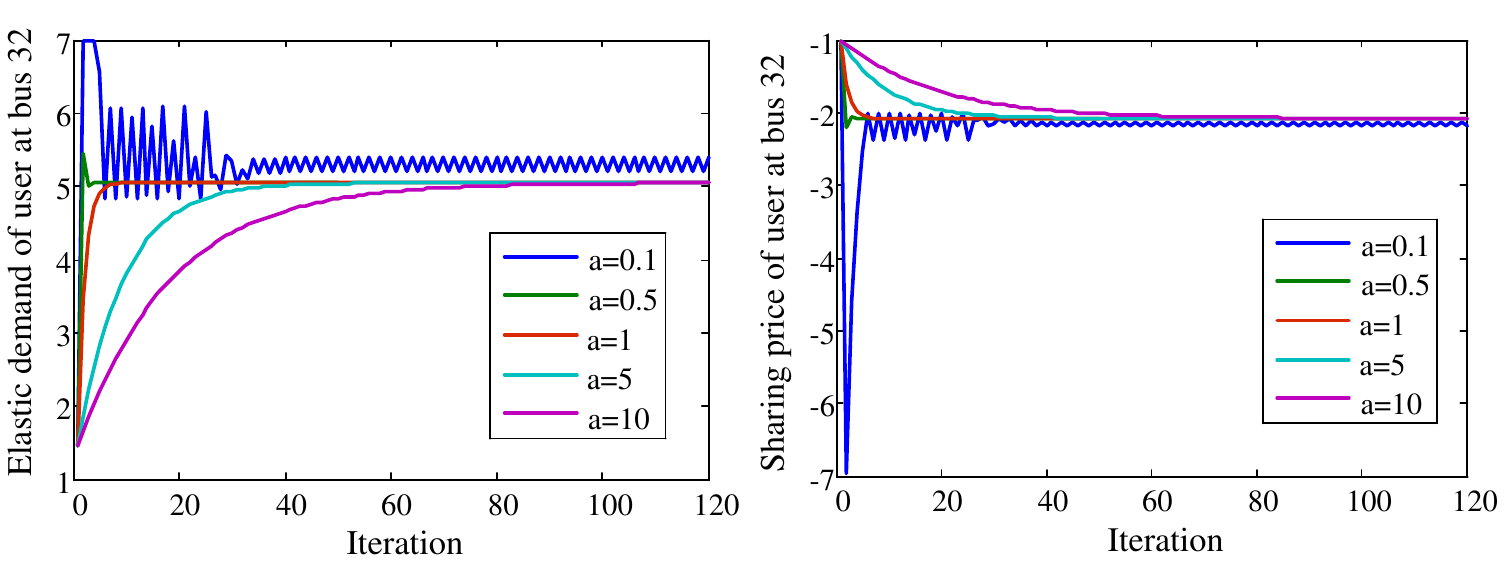}
	\caption{Change of strategy and price under different $a$.}
	\label{fig:impact_of_a}
\end{figure}

First, we test the impact of $a$ on the performance of the modified BR algorithm by the 69-bus system. Here, Condition C1 is satisfied when $a>3.84$. Let $a$ equals to 0.1, 0.5, 1, 5, and 10, respectively, and the changes of user's demand and the sharing price over iterations are plotted in Fig. \ref{fig:impact_of_a} (take the user at bus 32 as an example). We can find that, with a smaller $a$, the modified best response algorithm converges faster; however, when $a$ is too small, oscillation may occur. Moreover, when $a=0.5$ or 1, though Condition C1 is not met, the algorithm still converges, meaning that it is a sufficient but not necessary condition. Under all $a$, our algorithm reaches equilibrium within 180 seconds. 

We then test the performance of Algorithm 2 using the 38-bus and 69-bus cases with different settings in TABLE \ref{tab:performance}, where $F_l$ refers to the original flow limit. Generally, it takes a longer time and more iteration to output the absorbable region when the system constraints are more stringent. But under all cases, the algorithm terminates in less than a minute, which is acceptable for online use.
\begin{table}[h]
	\renewcommand{\arraystretch}{1.3}
	\centering
	\caption{Performance of Algorithm 2 under different settings.}
	\label{tab:performance}
	\begin{tabular}{cccc}
		\hline 
		Setting & Case 38 with $F_l$ & Case 38 with $2F_l$  & Case 38 with $4F_l$ \\
		Iteration & 31 & 20 & 15 \\
		Time (s) & 46.16 & 34.29 & 26.60\\
		\hline
		Setting & Case 69 with $F_l$ & Case 69 with $2F_l$  & Case 69 with $4F_l$ \\
		Iteration & 36 & 19 & 9 \\
		Time (s) & 49.69 & 27.27 & 19.17\\
		\hline
	\end{tabular}
\end{table}

\section{Conclusion}
With the mushrooming of distributed renewable energy at the demand side, new energy management schemes with explicit characterization on  flexibility are in great need. This paper proposes an energy sharing mechanism that encourages energy exchange among end-users. A generalized Nash game model is proposed to formulate the equilibrium of the sharing market. The energy transaction is coordinated via price signals which meet the requirement of information privacy. We prove that the generalized Nash equilibrium achieves social optimum, and leads to the same flexibility as the centralized dispatch in the sense of the absorbable region. The linear programming based projection algorithm can efficiently generate the boundaries of the absorbable region. Future research may use more accurate power flow models and consider multiple periods and temporal correlations.

\ifCLASSOPTIONcaptionsoff
\newpage
\fi
        
        \bibliographystyle{IEEEtran}
        \bibliography{IEEEabrv,mybib}

\appendix
\makeatletter
\@addtoreset{equation}{section}
\@addtoreset{theorem}{section}
\makeatother
\setcounter{equation}{0}  
\renewcommand{\theequation}{A.\arabic{equation}}
\subsection{Proof of Proposition \ref{prop-1}}
\label{apen-1}	
		Denote $p_k^{out}:=-a\lambda_k+b_k$, and let $\mathcal{Y}$ be the feasible region of $p^{out}_k,\forall k \in \mathcal{I} \cup \mathcal{J}$ characterized by \eqref{eq:sharing-oper.2} and \eqref{eq:sharing-oper.3} (or \eqref{eq:central.2} and \eqref{eq:central.5}). Obviously, $\mathcal{Y}$ is also a closed convex set.
		
		Then problem \eqref{eq:sharing-oper} can be equivalently written as
		\bsq
		\begin{align}
		\mathop{\min}_{p_k^{out},\forall k \in \mathcal{I} \cup \mathcal{J}}~ & \sum \limits_{k \in \mathcal{I} \cup \mathcal{J}} (p_k^{out}-b_k)^2 \\
		\mbox{s.t.}~ & \sum\limits_{k \in \mathcal{I} \cup \mathcal{J}} p_k^{out} =0 ; p^{out} \in \tilde{\mathcal{P}}
		\end{align}
		\esq
		Suppose $(d^*,b^*,\lambda^*)$ is a GNE of the game $\mathcal{G}(w)$, then for problem \eqref{eq:sharing-pro} we have
		\begin{align}
		\label{eq:condition-2}
		f_k(d_k)-f_k(d_k^*)+(d_k-d_k^*)\lambda_k^*\ge 0, \forall d \in \hat{\mathcal{D}}_k,\forall k \in \mathcal{I} \cup \mathcal{J}
		\end{align}
		Note that $b_k$ is eliminated by substituting constraint \eqref{eq:sharing-pro.2} into the objective function \eqref{eq:sharing-pro.1}. 
		For problem \eqref{eq:sharing-oper} we have
		\begin{align}
		\label{eq:condition-3}
		\sum \limits_{k \in \mathcal{I} \cup \mathcal{J}} (p_k^{out}-p_k^{out*})(p_k^{out*}-b_k^*) \ge 0,\forall p^{out} \in \mathcal{Y}
		\end{align}

		For problem \eqref{eq:central}, its Lagrangian function is
		\begin{align}
		L(d,p^{out},\lambda)=~&\sum \limits_{k \in \mathcal{I} \cup \mathcal{J}} f_k(d_k) +\sum \limits_{k \in \mathcal{I}} \lambda_k(d_k^f+d_k-p_k^{out}) \nonumber\\
		~& +\sum \limits_{k \in \mathcal{J}} \lambda_k(d_k^f+d_k-w_k-p_k^{out})
		\end{align}
		which is defined on $\Omega:=\prod_{k \in \mathcal{I} \cup \mathcal{J}} \hat{\mathcal{D}}_k \times \mathcal{Y} \times \mathbbm{R}^{(I+J)}$. Let $(\hat{d},\hat p^{out},\hat \lambda)$ be a saddle point of the Lagrangian function, then $(\hat d, \hat p^{out}, \hat \lambda) \in \Omega$ and it satisfies $\forall (d,p^{out},\lambda)\in \Omega$ \cite{kinderlehrer2000introduction}:
		\bsq
		\label{eq:condition}
		\begin{align}
		\left[f_k(d_k)- f_k(\hat d_k)+ (d_k-\hat d_k)\hat \lambda_k\right] & ~\ge 0 ,\forall k \in \mathcal{I} \cup \mathcal{J} \label{eq:condition.1}\\
		-\sum \limits_{k \in \mathcal{I} \cup \mathcal{J}} (p_k^{out}-\hat p_k^{out})(\hat \lambda_k) &~\ge 0 \label{eq:condition.2}\\
		\sum \limits_{k \in \mathcal{I}}  (\lambda_k-\hat \lambda_k)(d_k^f+ \hat d_k-\hat p_k^{out}) &\nonumber\\
		+\sum \limits_{k \in \mathcal{J}}  (\lambda_k-\hat \lambda_k)(d_k^f+ \hat d_k-w_k-\hat p_k^{out}) &~\le 0\label{eq:condition.3}
		\end{align}
		\esq
		
		\emph{\textbf{Existence}}. When Assumption A1 holds, suppose $\hat d$ is the optimal solution of \eqref{eq:central} and $\hat \lambda$ is the corresponding dual variable. Let $d^*=\hat d$, $\lambda^*=\hat \lambda$, $p_k^{out*}=d_k^f+\hat d_k$ for all $k$ in $\mathcal{I}$, $p_k^{out*}=d_k^f+\hat d_k-w_k$ for all $k$ in $\mathcal{J}$, and $b_k^*=a\hat \lambda_k+p^{out*}_k,\forall k \in \mathcal{I} \cup \mathcal{J}$, then it is easy to check that \eqref{eq:condition-2} and \eqref{eq:condition-3} are met. Thus, we have constructed a GNE $(d^*,b^*,\lambda^*)$.
		
		\emph{\textbf{Uniqueness}}. Given a GNE $(d^*,b^*,\lambda^*)$, when $k \in \mathcal{I}$, we have $b_k^*=d_k^f+d_k^*+a\lambda_k^*$; when $k \in \mathcal{J}$, we have $b_k^*=d_k^f+d_k^*-w_k+a\lambda_k^*$. Let $\hat d=d^*$, $\hat \lambda=\lambda^*$, and $\hat p^{out}=-a\lambda^*+b^*$, then it is easy to check that $(\hat d,\hat p^{out},\hat \lambda)$ satisfies \eqref{eq:condition}, so $\hat d$ is the optimal solution of \eqref{eq:central} and $\hat \lambda$ is the corresponding dual variable. Since the objective function is strictly convex, and the constraint sets $\hat{D}_k,\forall k \in \mathcal{I} \cup \mathcal{J}$  and $\tilde{\mathcal{P}}$ are all closed convex sets, problem \eqref{eq:central} has a unique solution \cite{LectureNote}, so $\hat d$ is unique.

\setcounter{equation}{0}  
\renewcommand{\theequation}{B.\arabic{equation}}
\subsection{Convergence of Modified Best-Response based Algorithm}
\label{apen-2}
Let $q_k:=-a\lambda_k+b_k$ for all $k \in \mathcal{I} \cup \mathcal{J}$. At the $n$-th iteration, given $b^{n}$, the microgrid operator's problem is equivalent to
\bsq
\label{eq:platform-update}
\begin{align}
q^{n+1}= ~&\mbox{argmin}\{\theta_1(q)-\frac{(b^{n})^Tq}{a} \nonumber\\
& +\frac{|q-q^n+b^{n-1}-b^{n}|^2}{2a} | ~q \in \mathcal{Q}\} \label{eq:platform-update.1}\\
\lambda^{n+1}= ~ & \frac{1}{a}(b^{n}-q^{n+1}) \label{eq:platform-update.2}
\end{align}
\esq
where $\theta_1(q):=\frac{1}{2a}\sum \nolimits_{k \in \mathcal{I} \cup \mathcal{J}} q_k^2$ and 
$$\mathcal{Q}:=\left\{q~|~\mbox{s.t.}~\sum \nolimits_{k \in \mathcal{I} \cup \mathcal{J}} q_k = 0, q \in \tilde{\mathcal{P}}\right\}$$.
The users’ problems \eqref{eq:sharing-pro} are equivalent to:
\bsq
\label{eq:prosumer-update}
\begin{align}
d^{n+1}=~&\mbox{argmin}\{\theta_2(d)+\frac{(b^n)^Td}{a}   \nonumber\\
~ & +\frac{|q^{n+1}-d-D|^2}{2a}| d \in \cup_{k \in \mathcal{I} \cup \mathcal{J}} \hat{\mathcal{D}}_k\} \label{eq:prosumer-update.1}\\
b^{n+1}=~&b^n-q^{n+1}+d^{n+1}+D \label{eq:prosumer-update.2}
\end{align}
\esq
where $\theta_2(p):=\sum \limits_{k \in \mathcal{I} \cup \mathcal{J}} f_k(d_k) - \frac{1}{2a}\sum \limits_{k \in \mathcal{I} \cup \mathcal{J}} (d_k+D_k)^2$, $D_k=d_k^f,\forall k \in \mathcal{I}$, $D_k=d_k^f,\forall k \in \mathcal{J}$.

If Condition C1 holds, both $\theta_1(q)$ and $\theta_2(d)$ are convex, and the following function has a unique saddle point $(q^*,d^*,b^*)$.
\begin{align}
\label{eq:Lfuntion}
\theta_1(q)+\theta_2(d)-\frac{b^T(q-d-D)}{a}
\end{align}
By variational inequality technique \cite{patriksson2013nonlinear} together with \eqref{eq:platform-update.2} and \eqref{eq:prosumer-update.2}, \eqref{eq:platform-update.1} is equivalent to for all $q \in \mathcal{Q}$:
\begin{align}
\label{eq:variational-eq1}
	\theta_1(q)-\theta_1(q^{n+1})+(q-q^{n+1})^T \left\{-\frac{1}{a}b^{n+1}+\frac{1}{a}(d^{n+1}-d^{n})\right\} \ge 0 
\end{align}  
Similarly, user’s problem is equivalent to for all $d \in \cup_{k \in \mathcal{I} \cup \mathcal{J}} \hat{\mathcal{D}}_k$:
\begin{align}
\label{eq:variational-eq2}
\theta_2(d)-\theta_2(d^{n+1})+(d-d^{n+1})^T \left\{\frac{1}{a}b^{n}-\frac{1}{a}(q^{n+1}-d^{n+1}-D)\right\} \ge 0 
\end{align} 
Let $t:=(q,d)$ and $\theta(t):=\theta_1(q)+\theta_2(d)$, \eqref{eq:variational-eq1}-\eqref{eq:variational-eq2} implies
\begin{align} 
\theta(t)-\theta(t^{n+1})+ 
\left(                
\begin{array}{c}   
q-q^{n+1} \\ 
d-d^{n+1} \\
b-b^{n+1} \\
\end{array}
\right)^T \cdot ~& \nonumber\\
\left\{\left(                
\begin{array}{c}   
-b^{n+1}/a \\ 
b^{n+1}/a \\
q^{n+1}\!-\!d^{n+1}\!-\!D \\
\end{array}
\right)+
\left( 
\begin{array}{c}   
-\mathbf{I}/a \\
\mathbf{I}/a \\
0\\
\end{array}
\right)(d^{n}-d^{n+1})
\right.  ~& \nonumber\\
\left.+
\left(                
\begin{array}{cc}   
0 & 0 \\
{\mathbf{I}/a} & 0 \\
0 & \mathbf{I}/a\\
\end{array}
\right)
\left(                
\begin{array}{c}   
d^{n+1}-d^{n} \\
b^{n+1}-b^{n} \\
\end{array}
\right)\right\}  ~& \ge 0 \nonumber
\end{align}
Let $w:=(q,d,b) \in \mathcal{Q} \times \cup_{k \in \mathcal{I} \cup \mathcal{J}} \hat{\mathcal{D}}_k \times \mathbb{R}^{|\mathcal{I}|\times |\mathcal{J}|}$ and define a mapping $F(w):=(-b/a,b/a,q-d-D)$, which is indeed monotone. Then we have
\begin{align}
& \theta(t^{n+1})-\theta(t^*)+(w^{n+1}-w^*)^TF(w^{n+1}) \nonumber\\
\ge ~&\theta(t^{n+1})-\theta(t^*)+(w^{n+1}-w^*)^TF(w^*)  { \geq 0} \nonumber
\end{align}
Therefore, we have
\begin{align}
\label{eq:10}
&   \left( 
\begin{array}{c}   
d^*- d^{n+1} \\
b^*- b^{n+1} \\
\end{array}
\right)^T
\left( 
\begin{array}{cc}   
\mathbf{I}/a & 0 \\
0 &  \mathbf{I}/a \\
\end{array}
\right)
\left( 
\begin{array}{c}   
d^{n+1}-d^{n} \\
b^{n+1}-b^{n} \\
\end{array}
\right) \nonumber\\
\ge ~ & (w^{n+1}-w^*)^T \left( 
\begin{array}{c}   
-\mathbf{I}/a \\
\mathbf{I}/a \\
0\\
\end{array}
\right)(d^{n}-d^{n+1}) \nonumber \\
=~& \frac{1}{a}(b^{n+1}-b^{n})^T(d^n-d^{n+1}) \ge 0
\end{align}
Note that the last inequality is because of \eqref{eq:variational-eq2}. We could get
\begin{align}
\label{eq:converge}
\bigg| 
\begin{array}{c}   
d^n-d^{*} \\
b^n-b^{*} \\
\end{array}
\bigg|^2 \ge~& 
\bigg| 
\begin{array}{c}   
d^{n+1}-d^{*} \\
b^{n+1}-b^{*} \\
\end{array}
\bigg|^2 + 
\bigg| 
\begin{array}{c}   
d^n-d^{n+1} \\
b^n-b^{n+1} \\
\end{array}
\bigg|^2 
\end{align}
The sequence $\{(d^n, b^n)\}$ is F\'ejer monontone \cite{combettes2001fejer}, with $|(d^n-d^*)^T, (b^n-b^*)^T|^2$ decreasing in each iteration $n$ by $|(d^n-d^{n+1})^T, (b^n-b^{n+1})^T|^2$. As a result, the sequence $\{|(d^n-d^*)^T, (b^n-b^*)^T|^2\}$ converges and sequences $\{d^n\}$ and $\{b^n\}$ are bounded. With \eqref{eq:converge}, the sequence $\{d^n\}$ ($\{b^n\}$) only has one cluster point. According to \eqref{eq:variational-eq2} we can get $d^n \to d^*$ and $b^n\to b^*$. With \eqref{eq:platform-update} we know $q^n \to q^*$ and $\lambda^n \to \lambda^*$

\end{document}